\documentclass[12pt]{article}
\usepackage[centertags]{amsmath}
\usepackage{amsfonts}
\usepackage{amssymb}
\usepackage{latexsym}
\usepackage{amsthm}
\usepackage{newlfont}
\usepackage{graphicx}
\usepackage{listings}
\usepackage{booktabs}
\usepackage{abstract}
\usepackage{enumerate}
\usepackage{xcolor}
\usepackage[margin=1in]{geometry}
\date{}
\usepackage{fancyhdr}
\usepackage[hidelinks]{hyperref}
\hypersetup{pdftitle={Basis selection, universal sampling discretization, and sampling recovery},
pdfauthor={Feng Dai, Jiashu Ran, Vladimir N. Temlyakov}}
\usepackage{multirow} 
\newlength{\defbaselineskip}
\newcommand{\setlinespacing}[1]%
           {\setlength{\baselineskip}{#1 \defbaselineskip}}

{\medskip\noindent \textit{Proof of #1. }}%
{\qed \medskip}

\def\cD{{\mathcal D}}

\def\cG{{\mathcal G}}

\def\cM{{\mathcal M}}

\def\cX{{\mathcal X}}

\def\bbC{{\mathbb C}}

\def\bbN{{\mathbb N}}

\def\bbR{{\mathbb R}}

\def\bbZ{{\mathbb Z}}

\def\bA{{\mathbf A}}

\def\bF{{\mathbf F}}

\def\bW{{\mathbf W}}

\def\bk{\mathbf k}

\def\bs{\mathbf s}
\def\bt{\mathbf t}

\def\bu{\mathbf u}
\def\bv{\mathbf v}

\def\bx{\mathbf x}
\def\by{\mathbf y}

 \def \<{\langle}
\def\>{\rangle}

\def \Og{\Omega}

\def \ta{\theta}
\def\th{\theta}

\def\bt{\beta}

\def\la{\lambda} 

\def\ld{\lambda}
\def\RR{\mathbb{R}}

\def\sa{\sigma}
\def\SS{\mathbb{S}}
\def\EE{\mathbb{E}}
\def\PP{\mathbb{P}}

\def\vi{\varphi}

\def\T{\mathsf T}

\def\f{\frac}

\def\bt{\beta}

\def\CD{\mathcal{D}}
\def\Bl{\Bigl}
\def\Br{\Bigr}
\def\wt{\widetilde}
\newtheorem{Theorem}{Theorem}[section]
\newtheorem{Lemma}[Theorem]{Lemma}
\newtheorem{Definition}[Theorem]{Definition}
\newtheorem{Proposition}[Theorem]{Proposition}
\newtheorem{Remark}[Theorem]{Remark}

\newtheorem{Corollary}[Theorem]{Corollary}
\numberwithin{equation}{section}

\newcommand{\be}{\begin{equation}}
\newcommand{\ee}{\end{equation}}

\def\eps{\epsilon}

\def\la{\langle}
\def\ra{\rangle}

\newcommand{\eq}[1]{\begin{align*}#1\end{align*}}

\newcommand{\beq}{\begin{equation}}
\newcommand{\enq}{\end{equation}}

\begin{document}

\title{Basis selection, universal sampling discretization, and sampling recovery}

\author{Feng Dai, Jiashu Ran, and Vladimir N. Temlyakov}

\newcommand{\Addresses}{{
  \bigskip
  \footnotesize

  Feng Dai, \textsc{ Department of Mathematical and Statistical Sciences\\
University of Alberta\\ Edmonton, Alberta T6G 2G1, Canada\\
E-mail:} \texttt{fdai@ualberta.ca }

\medskip

 Jiashu Ran, \textsc{ Department of Mathematical and Statistical Sciences\\
University of Alberta\\ Edmonton, Alberta T6G 2G1, Canada\\
E-mail:} \texttt{jran@ualberta.ca }

 \medskip
  Vladimir N. Temlyakov, \textsc{University of South Carolina,\\ Steklov Institute of Mathematics,\\  Lomonosov Moscow State University,\\ and Saint Petersburg State University.
  \\
E-mail:} \texttt{temlyak@math.sc.edu}

}}
\maketitle

\begin{abstract}
For every real $N$-dimensional subspace of $L_2$, we construct one
orthonormal basis such that, for every $1\le v\le N$, at most
$Cv\log(eN/v)$ positively weighted samples discretize the squared
$L_2$ norm of all $v$-term combinations. Random rotations also give
simultaneous weighted $L_p$ estimates and a sharp $C\log(eN)$ bound
for the squared $L_2$ norm of the pointwise maximum. For prescribed
coefficient classes, a change of measure yields sampling recovery
estimates from $Cv\log^3(2v)$ samples under a lower Riesz inequality
and a common square-integrable pointwise bound.
\end{abstract}



\section{Introduction}\label{In}

We study how the choice of a basis and a sampling measure affects
universal discretization and sampling recovery.  Our first result concerns
an arbitrary real $N$-dimensional subspace of $L_2$. We construct one
orthonormal basis for which, at every sparsity level $1\le v\le N$,
the squared $L_2$ norm of all $v$-term combinations can be discretized
using at most $Cv\log(eN/v)$ positively weighted samples. The basis is
independent of $v$.

To be precise, let $(\Omega,\mu)$ be a probability space. We write 
\[
\|f\|_p:=\|f\|_{L_p(\mu)}=\left( \int_{\Omega} |f|^p \, d\mu \right)^{\frac{1}{p}}, \qquad 1 \leq p <\infty.
\]
For a finite system $\mathcal{D}_N=\{g_j\}_{j=1}^{N}$, put 
\[
\Sigma_v (\mathcal{D}_N): = \left\{ \sum_{j\in \Lambda} c_j g_j:\, \Lambda \subset \{1,\ldots, N\}, \, |\Lambda| \leq v \right\}, \quad 1\leq v \leq N.
\]  
Thus, $\Sigma_v(\cD_N)$ is the union of the subspaces generated by
at most $v$ elements of $\cD_N$. Universal sampling discretization
requires a single set of sampling points and weights that
discretizes the norm simultaneously on all these subspaces.

Let $X_N$ be a real $N$-dimensional subspace of $L_2(\mu)$ with $N\ge 2$.
To make pointwise evaluation well defined, we fix a basis of $X_N$
and choose finite measurable representatives of its elements
defined everywhere on $\Omega$, modifying them on null sets
if necessary. Every $f\in X_N$ is then represented pointwise
by the corresponding linear combination of these basis functions.
We use these representatives throughout. Continuity is not assumed.
Define the Christoffel function by
\begin{equation}\label{2-1}
 \lambda_N(\bx):=
 \inf_{\substack{f\in X_N\\f(\bx)=1}}\|f\|_2^2,
 \qquad \bx\in\Omega,
\end{equation}
with the convention that $\inf\emptyset=\infty$.

Our main discretization result is the following.

\begin{Theorem} \label{thm-1-1}Let  $X_N$ be an $N$-dimensional space of bounded functions on $\Og$, with $N\ge 2$,  satisfying
		\begin{equation*}
			\sup_{\bx\in\Og}|f(\bx)| \leq \sqrt{KN}  \|f\|_{L_2(\mu)},\   \ \forall f\in X_N
		\end{equation*}
for some constant $K\ge 1$. 
{Then for every $0<\epsilon<1$, there exists} an orthonormal basis
$\CD_N=\{\vi_j\}_{j=1}^N$ of $(X_N,\|\cdot\|_{L_2(\mu)})$ such that,
for every integer $1\le v\le N$, there are points
$\xi_v^1,\ldots,\xi_v^m\in{\Omega}$, {not necessarily distinct,} satisfying
\begin{equation}\label{UDRunweightedCount}
 m\le CK\epsilon^{-2}v\log\Bl(\f {eN}v\Br)
\end{equation}
and
\begin{equation}\label{UDRunweightedDiscretization}
 (1-\epsilon)\|f\|_{L_2(\mu)}^2
 \le\f 1 m \sum_{j=1}^m |f(\xi_v^j)|^2
 \le(1+\epsilon)\|f\|_{L_2(\mu)}^2,
 \qquad f\in\Sigma_v(\CD_N).
\end{equation}
Here $C$ is an absolute constant, and for each fixed $\epsilon$, the basis
is independent of $v$.
\end{Theorem}

\begingroup

Of particular interest is the case $v=N$, in which
$\Sigma_N(\CD_N)=X_N$ is independent of the choice of basis.
Limonova and Temlyakov \cite[Theorem~1.1]{LT2022} proved a
discretization inequality of the form
\eqref{UDRunweightedDiscretization}, with $m\le CKN$ and with
$1-\epsilon$ and $1+\epsilon$ replaced by $c_1$ and $c_2K$,
respectively, where $C,c_1,c_2>0$ are absolute constants.
Their proof uses an iteration argument based on the theorem of
Marcus, Spielman, and Srivastava \cite{MSS2015}.
Refining this argument, Kosov \cite[Theorem~2.4]{Kosov2022}
obtained the lower and upper constants $1-\epsilon$ and
$K(1+\epsilon)$ for every $0<\epsilon<1$, with
$m\le C\epsilon^{-2}KN$.

For general $v$, our proof uses Bownik's selection theorem
\cite[Theorem~1.5]{B26} to construct sampling matrices (see \eqref{Sampling_Matrix}) with
estimates for both their operator norms and traces.  We then
apply spherical concentration after a random orthogonal rotation
to obtain one basis that works simultaneously for all
$1\le v\le N$.
\par\endgroup

{A standard change of density gives the following
consequence of Theorem~\ref{thm-1-1}.}
\begin{Corollary}\label{UDRchosen}
 Let $X_N$ be an $N$-dimensional subspace of $L_2(\mu)$ with $N\ge2$, and let  
  $0<\epsilon<1$. Then there exists an orthonormal basis
$\CD_N=\{\vi_j\}_{j=1}^N$ of $(X_N,\|\cdot\|_{L_2(\mu)})$ such that,
for every integer $1\le v\le N$, there are points
$\xi_v^1,\ldots,\xi_v^m\in\Omega$  satisfying
\begin{equation}\label{UDRchosenCount}
 m\le C\epsilon^{-2}v\log\Bl(\f{eN}v\Br)
\end{equation}
and
\begin{equation}\label{UDRdiscretization}
 (1-\epsilon)\|f\|_{L_2(\mu)}^2
 \le\f N m \sum_{j=1}^m \ld_N(\xi_v^j)|f(\xi_v^j)|^2
 \le(1+\epsilon)\|f\|_{L_2(\mu)}^2,
 \qquad f\in\Sigma_v(\CD_N),
\end{equation}
where $C$ is an absolute constant, and it is agreed that $0\cdot\infty=0$.
\end{Corollary}

\begingroup

Corollary~\ref{UDRchosen} requires no boundedness assumption on
the original basis of $X_N$. The case $v=N$ corresponds to the
weighted discretization result of Kosov
\cite[Corollary~2.5]{Kosov2022}, which gives
\eqref{UDRdiscretization} with prescribed accuracy
$\epsilon\in(0,1)$ and $m\le C\epsilon^{-2}N$.
For each fixed $\epsilon$, Corollary~\ref{UDRchosen} provides
one orthonormal basis for which the discretization inequality
holds at every sparsity level $1\le v\le N$, with
$m\le C\epsilon^{-2}v\log(eN/v)$. The sampling points and
weights may depend on $v$.
\par\endgroup

\begingroup
The random orthogonal rotations used to obtain the discretization result in Theorem~\ref{thm-1-1} also
allow us to control the weighted $L_p$ estimates of the basis. This connection
is developed in Section~\ref{OB}. More precisely, Theorem~\ref{OBT1} shows that,
after normalization by the Christoffel function, a random orthogonal
rotation yields simultaneous $L_p$ bounds of order $\sqrt p$ for all $p>2$.  These bounds can be imposed on a basis satisfying
Corollary~\ref{UDRchosen}; see Corollary~\ref{UDRallvLq}.
The same basis can also be chosen so that the integral of its squared
pointwise maximum is bounded by $C\log(eN)$; see Corollary~\ref{UDRlogenv}. This logarithmic order is sharp in general,
as shown by the space of spherical harmonics of degree one;
see Proposition~\ref{UDRsharpEnvelope}.

\par\endgroup

\begingroup
However, a change of basis can alter the coefficient conditions defining a
function class.  We therefore treat recovery for prescribed coefficient
classes separately. In Section~\ref{UD}, we use a change of measure
to obtain a uniformly bounded system while preserving the expansion
coefficients and sampling recovery errors. Universal discretization
of the transformed system then yields weighted discretization
inequalities for the original one, which provide the
foundation for the recovery estimates in the following Section \ref{SR}.
\par\endgroup

\begingroup
For these results in Section~\ref{SR}, we start with some necessary definitions and notation. Given a function class $\bF$ with specified point values, define
\[
 \varrho_m^o(\bF,L_p)
 :=\inf_{\xi\in\Omega^m}\inf_{\cM}
       \sup_{f\in\bF}
       \|f-\cM(f(\xi^1),\ldots,f(\xi^m))\|_p,
\]
where $\cM:\bbC^m\to L_p(\mu)$ is an arbitrary reconstruction map; for
real classes one may use $\bbR^m$. Repetitions in the sampling tuple are
allowed. Thus the definition also permits at most $m$ distinct sampling
points. The superscript $o$ refers to optimal recovery in the sense of
information-based complexity \cite{TWW}.  

\par\endgroup

\begingroup
Let $\cD=\{g_\bk\}_{\bk\in\bbZ^d}$, and let $\cG=\{G_j\}_{j\ge0}$ be
finite index sets satisfying
\be\label{Di4}
 G_0=\{\mathbf0\}\subset G_1\subset G_2\subset\cdots,
 \qquad \bigcup_{j\ge0}G_j=\bbZ^d.
\ee
For $0<\bt\le1$, $r>0$, and $b\in\bbR$, denote by
$\bA^{r,b}_\bt(\cD,\cG)$ the class of functions admitting a representation
\be\label{Direpr}
 f=\sum_{\bk\in\bbZ^d}a_\bk g_\bk,
 \qquad \sum_\bk|a_\bk|<\infty,
\ee
such that
\be\label{DiAr}
 |a_{\mathbf0}|\le1,
 \qquad
 \left(\sum_{\bk\in G_j\setminus G_{j-1}}|a_\bk|^\bt\right)^{1/\bt}
 \le2^{-rj}j^b,\qquad j\ge1.
\ee
These conditions describe decay from one block to the next and sparsity
within each block. They are imposed on an admissible representation;
uniqueness is not required. If a coefficient quasi-norm is used for a
redundant system, it is the infimum over all admissible representations.
The block $\ell_1$ norms are bounded by $2^{-rj}j^b$, whose sum is finite.
Hence a common finite pointwise bound gives absolute convergence
at every sampling point. If this bound belongs to $L_p$, the
series also converges in $L_p$.
Related Wiener--Beurling spaces are discussed in \cite{bull2020}.

\par\endgroup

We assume the regularity condition
\be\label{sa11}
 \mathrm{RC}(\th,\th'):\qquad
 2^{\th'}\le\frac{|G_{j+1}|}{|G_j|}\le2^\th,
 \quad j\ge1,\qquad 0<\th'<\th,
\ee
and put
\be\label{sa2}
 n(v):=\min\{j\ge1:\ v\le|G_j|\},\qquad v\ge1.
\ee
Thus the cardinalities $|G_j|$ grow exponentially, and $n(v)$
specifies the level at which the initial index set first contains at
least $v$ elements.

Typical examples are the dyadic cubes
\[
 G_j^c:=\{\bk\in\bbZ^d:\|\bk\|_\infty<2^j\},\quad j\ge1,
 \qquad \bA_\bt^r(\cD):=\bA^{r,0}_\bt(\cD,\cG^c),
\]
and the dyadic hyperbolic crosses. For $\bs\in\bbN_0^d$, set
\[
 \Delta(\bs):=\{\bk\in\bbZ^d:
 \lfloor2^{s_i-1}\rfloor\le|k_i|<2^{s_i},\ 1\le i\le d\},
 \qquad G_j^h:=\bigcup_{\|\bs\|_1\le j}\Delta(\bs).
\]
The corresponding classes
$\bW^{a,b}_{A_\bt}(\cD):=\bA^{a,b}_\bt(\cD,\cG^h)$ are related to
classes with mixed smoothness, while the cube classes are related to
Sobolev--Nikol'skii--Besov classes. We refer to
\cite[Chapters~3--7]{VTbookMA} for the approximation background and to
\cite{STT,VT205,VT203,VT202} for sampling recovery on classes with
structural conditions.

We impose the following condition on the system $\Phi=\{\phi_{\mathbf{k}}\}_{\mathbf{k}\in \mathbb{Z}^d} \subset L_2(\mu)$. 

\begin{Definition}\label{InD1}
The system $\Phi$ satisfies a lower Riesz inequality with constant $K>0$
if, for every finite $\Gamma\subset\bbZ^d$ and every choice of coefficients,
\be\label{Bessel}
 \sum_{\bk\in \Gamma}|a_\bk|^2
 \le K\left\|\sum_{\bk\in \Gamma}a_\bk\phi_\bk\right\|_2^2.
\ee
\end{Definition}
This inequality is sometimes called a Bessel inequality in the sampling
literature. We use the term ``lower Riesz inequality'' to distinguish it
from the usual upper Bessel inequality.

Our recovery result applies to the structural classes $\mathbf{A}_{\beta}^{r,b}(\Phi,\mathcal{G})$.
\begin{Theorem}\label{InT1}
Let $\Phi=\{\phi_{\mathbf{k}}\}_{\mathbf{k} \in \mathbb{Z}^{d}} \subset L_2 (\mu)$ be a system.
Let $r>0$, $b\in\bbR$, $0<\bt\le1$, and assume \eqref{Bessel} and  \eqref{sa11}. Suppose that there is a common measurable pointwise bound $B \in L_2(\mu)$ such that 
\[
|\phi_{\mathbf{k}}(\mathbf{x})| \leq B(\mathbf{x}), \quad  1 \leq B(\mathbf{x}) <\infty, \quad \mathbf{x} \in \Omega, \quad \mathbf{k} \in \mathbb{Z}^d. 
\]
There are constants
$c,C>0$, depending only on
$K,\|B\|_2,r,b,\bt,\th,\th'$ and $|G_1|$, such that, for every integer
$v\ge1$,
\be\label{In1a}
 \varrho_m^o(\bA^{r,b}_\bt(\Phi,\cG),L_2(\mu))
 \le C v^{1/2-1/\bt}2^{-rn(v)}n(v)^b
\ee
whenever $m\ge c v\log^3(2v)$. The reconstruction takes values in the span
of an initial dictionary containing at most $C(2v)^c$ elements. The same error bound holds
in $L_p(\mu)$ for $1\le p\le2$.
\end{Theorem}

No upper Riesz inequality is required. The proof in Section~\ref{SR} combines sparse approximation with universal $L_2$ discretization. The discretization estimate used here follows from a normalization of \cite[Theorem 1.1]{BDJR}; see Theorem~\ref{UDTsharp}. If $2<p<\infty$ and the same function $B$ belongs to $L_p(\mu)$, then Theorem~\ref{SRTp} gives a corresponding $L_p$ estimate, with the factor $v^{1/2-1/\beta}$ replaced by $v^{1-1/p-1/\beta}$. For the dyadic cube classes, Corollary~\ref{SRT7} expresses the error directly in terms of the number of samples. These are upper estimates; optimality of the logarithmic factors is not asserted.

For completeness, Appendix~\ref{Lp} records the corresponding norming-functional identities and the classical Auerbach construction.

The paper is organized as follows. Section~\ref{UDR} proves the
discretization theorem with one basis for all sparsity levels.
Section~\ref{OB} establishes simultaneous $L_q$ bounds for random
orthonormal bases and a sharp logarithmic bound for the integral
of their squared pointwise maximum. Section~\ref{UD} develops a
change of measure and weighted universal discretization, including
independent sampling from the resulting density. Section~\ref{SR} proves the recovery
estimates. Appendix~\ref{Lp} treats dual orthogonal systems.

\section{Universal discretization with a suitable basis}\label{UDR}

\begingroup

Our main purpose in this section is to prove Theorem~\ref{thm-1-1} and Corollary~\ref{UDRchosen}.
\par\endgroup
\begingroup

We work with a real $N$-dimensional subspace $X_N$ of $L_2(\mu)$,
where $N\ge2$. Corollary~\ref{UDRchosen} also holds for complex
coefficient vectors in the complex span of the real basis.
\par\endgroup
 
\subsection{Preliminary lemmas}
\begingroup

Let $\SS^{n-1}:=\{\bx\in\RR^n:|\bx|=1\}$ be the unit sphere in
$\RR^n$, equipped with the normalized surface measure $\sigma_n$.
Here and throughout this section, $|\cdot|$ denotes the Euclidean
norm, and $\|\cdot\|$ denotes the corresponding operator norm.
We first recall the following spherical concentration inequality.
\par\endgroup

\begin{Lemma}[{\cite[Theorem~5.1.4]{VR}}]\label{concentration}
Let $n\ge2$, and suppose that $f:\SS^{n-1}\to\RR$ satisfies
\[
 |f(\bx)-f(\by)|\le L|\bx-\by|,
 \qquad \bx,\by\in\SS^{n-1},
\]
for some $L>0$. There is an absolute constant $C>0$ such that
\begin{align}
 \sigma_n\left\{\bx:\left|f(\bx)-\int_{\mathbb{S}^{n-1}} f\,d\sigma_n\right|
                  >C\tau\right\}
 &\le2\exp\left(-\frac{n\tau^2}{2{L}^2}\right),
 \qquad \tau>0.
\end{align}
\end{Lemma}

Let $\mu_n$ denote the normalized Haar measure on the rotation group
$\mathrm{SO}(n)$. This measure is invariant under left and right
multiplication by a fixed rotation, and under transposition. 
The following identity, which relates Haar measure to spherical measure, is well known.

\begin{Lemma}\label{HaarSphere}
Let $n\ge2$ and $\mathbf e\in\SS^{n-1}$. Then, for every
$f\in L_1(\SS^{n-1},\sigma_n)$,
\begin{equation}\label{OB6}
 \int_{\mathrm{SO}(n)}f(Q\mathbf e)\,d\mu_n(Q)
 =\int_{\SS^{n-1}}f\,d\sigma_n.
\end{equation}
\end{Lemma}

\begingroup

For two real symmetric matrices $A$ and $B$ of the same size,
we write $A\preceq B$ if $B-A$ is positive semidefinite.
We also write $[m]:=\{1,\ldots,m\}$ for a positive integer $m$.
\par\endgroup

\begingroup

We need the following finite-dimensional form of Bownik's theorem.
\begin{Lemma}[{\cite[Theorem~1.5]{B26}}]\label{lem-2-3}
Let $\eps>0$, and let $T_1,\ldots,T_{2r}$ be positive semidefinite
$N\times N$ matrices satisfying
\[
 T:=\sum_{i=1}^{2r}T_i\preceq I_N\  \ \text{and}\  
 \ \operatorname{tr}T_i\le\eps,
 \qquad \text{for all $i$}.
\]
Then for every partition $\{J_k\}_{k=1}^r$ of $[2r]$ into two-element
sets, there exists a set $J\subset[2r]$ such that
$|J\cap J_k|=1$ for every $k$, and
\beq
 \left\|\sum_{i\in J}T_i-\f12T\right\|
 =\left\|\sum_{i\in J^c}T_i-\f12T\right\|
 \le2\sqrt{\eps}+\eps,
\enq
where $J^c=[2r]\setminus J$.
\end{Lemma}

As a direct consequence of Lemma~\ref{lem-2-3}, we obtain
the following selection lemma.
\begin{Lemma}\label{lem-2-4}
Let $u_1,\ldots,u_{2r}\in\RR^N$ satisfy $|u_i|^2\le R$
for all $i$, where $R>0$, and put
$
 S:=\f1{2r}\sum_{i=1}^{2r}u_iu_i^{\mathsf T}.
$
Then there exists a set $J\subset[2r]$ with $|J|=r$ such that
\[
 S_J:=\f1{|J|}\sum_{i\in J}u_iu_i^{\mathsf T}
\]
satisfies
\beq\label{2-4-0}
 \|S_J-S\|\le2\sqrt{2}\sqrt{\f{R\|S\|}{r}}+\f Rr
\enq
and
\beq\label{2-4}
 |\operatorname{tr}S_J-\operatorname{tr}S|\le\f R{2r}.
\enq
\end{Lemma}

\begin{proof}
Without loss of generality, we may
assume that $S\ne0$, and
\[
 |u_1|^2\le|u_2|^2\le\cdots\le|u_{2r}|^2.
\]
Set $T_i:=u_iu_i^{\mathsf T}/(2r\|S\|)$ for $1\le i\le2r$.
Then
\[
 \sum_{i=1}^{2r}T_i=\f S{\|S\|}\preceq I_N\  \ \text{and}\  \  \operatorname{tr}T_i
 =\f{|u_i|^2}{2r\|S\|}\le\f R{2r\|S\|}.
\]
Using Lemma~\ref{lem-2-3} with $\eps=R/(2r\|S\|)$, we obtain a set
$J\subset[2r]$ such that $|J\cap \{2j-1,2j\}|=1$  for every integer $1\le j\le r$, and 
\[
 \left\|\sum_{i\in J}\f{u_iu_i^{\mathsf T}}{2r\|S\|}
       -\f12\f S{\|S\|}\right\|
 \le2\sqrt{\f R{2r\|S\|}}+\f R{2r\|S\|}.
\]
Multiplying by $2\|S\|$ gives \eqref{2-4-0}.

To prove \eqref{2-4}, write
\[
 \{2j-1,2j\}=\{i_j,i_j'\},\  \ \text{where}\  \ i_j\in J,  \  \  i_j'\in J^c,
 \ \  1\le j\le r.
\]
Then
\begin{align*}
 \operatorname{tr}(S_J-S)
 &= \f1{r}\sum_{i\in J}|u_i|^2
                    - \f1{2r}\sum_{i=1}^{2r}|u_i|^2=\f1{2r}\left(\sum_{i\in J}|u_i|^2
                    -\sum_{i\in J^c}|u_i|^2\right)\\
                    &=\f1{2r}\sum_{j=1}^r\bigl(|u_{i_j}|^2-|u_{i_j'}|^2\bigr).
\end{align*}
It  follows that
\begin{align*}
 |\operatorname{tr}(S_J-S)|
 &\le\f1{2r}\sum_{j=1}^r\bigl(|u_{2j}|^2-|u_{2j-1}|^2\bigr)
 \le\f1{2r}\sum_{i=2}^{2r}\bigl(|u_i|^2-|u_{i-1}|^2\bigr)\\
 &=\f1{2r}\bigl(|u_{2r}|^2-|u_1|^2\bigr)\le\f R{2r}.
\end{align*}
This proves the lemma.
\end{proof}

We now apply Lemma~\ref{lem-2-4} repeatedly to construct
a sampling matrix with a prescribed number of points.
\begin{Lemma}\label{lem-sampling-matrices}
Let $(\Omega,\mu)$ be a probability space, and let
$\Phi:\Omega\to\RR^N$ be bounded and measurable. Assume that
\[
 \int_\Omega\Phi(\bx)\Phi(\bx)^{\mathsf T}\,d\mu(\bx)=I_N\  \ \text{and}\  \  |\Phi(\bx)|^2\le KN,
 \qquad \forall \bx\in\Omega,
\]
where $K\ge1$. Then for every integer $m\ge1$, there are points
$\bx_1,\ldots,\bx_m\in\Omega$, with repetitions allowed, such that
\begin{equation}\label{Sampling_Matrix}
    A:=\f1m\sum_{j=1}^m\Phi(\bx_j)\Phi(\bx_j)^{\mathsf T}
\end{equation}
satisfies
 \beq\label{2-6a}
 \left|\operatorname{tr}(A-I_N)\right|\le\f{CKN}{m},
\enq
and
\beq\label{2-7a}
 \|A-I_N\|\le C\left(\sqrt{\f{KN}{m}}+\f{KN}{m}\right),
\enq
where $C$ is an absolute constant.
\end{Lemma}

\begin{proof}Set $R:=KN$.  By a general
$L_2$-discretization result (see \cite[Theorem~3.2]{DKT2026}), if $s$ is an integer satisfying  $s\ge CRm\log N$ for a
sufficiently large absolute constant $C$, then 
there exist $M:=2^sm$ points $\by_1,\ldots,\by_M\in\Omega$
such that \beq
 A_0:=\f1M\sum_{j=1}^M\Phi(\by_j)\Phi(\by_j)^{\mathsf T}
 \quad\text{satisfies}\quad
 \|A_0-I_N\|\le\f1{2m}.
\enq
In the sequel, we fix such an integer $s$ and points $\by_1,\ldots, \by_M$.

Set   $J_0=[M]$. Applying Lemma~\ref{lem-2-4} successively, we obtain
index sets
\[
 J_s\subset J_{s-1}\subset\cdots\subset J_0,\  \ \text{with}\  \  |J_i|=2^{s-i}m,
\]
such that  the matrices
\[
 A_i:=\f1{|J_i|}\sum_{j\in J_i}
       \Phi(\by_j)\Phi(\by_j)^{\mathsf T},
 \qquad 0\le i\le s,
\]
satisfy, for $1\le i\le s$,
\beq\label{2-5}
 \|A_i-A_{i-1}\|
 \le2\sqrt{2}\sqrt{\f{R\|A_{i-1}\|}{2^{s-i}m}}
       +\f R{2^{s-i}m}
\enq
and
\beq\label{2-6}
 |\operatorname{tr}(A_i-A_{i-1})|
 \le\f R{2^{s-i+1}m}.
\enq

For $0\le i\le s$, put $e_i:=\|A_i-I_N\|$.
Since $\|A_{i-1}\|\le1+e_{i-1}$, \eqref{2-5} gives
\begin{align*}
 e_i
 &\le \|A_i-A_{i-1}\|+e_{i-1}\leq  e_{i-1}+2\sqrt{\f{2R}{2^{s-i}m}}\cdot\sqrt{1+e_{i-1}}
                   +\f R{2^{s-i}m}\\
 &\le\left(\sqrt{1+e_{i-1}}+
                  \sqrt{\f{2R}{2^{s-i}m}}\right)^2-1.
\end{align*}
Consequently,
\[
 \sqrt{1+e_i}-\sqrt{1+e_{i-1}}
 \le\sqrt{\f{2R}{2^{s-i}m}}.
\]
Summing these inequalities and using $e_0\le1/(2m)$, we obtain
\begin{align*}
 \sqrt{1+e_s}
 &\le\sqrt{1+e_0}+\sqrt{\f{2R}m}
                            \sum_{i=1}^s2^{-(s-i)/2}\\
 &\le1+\f1{4m}+2(\sqrt{2}+1)\sqrt{\f Rm}\le1+(2\sqrt{2}+3)\sqrt{\f Rm}.
\end{align*}
Squaring gives
\[
 \|A_s-I_N\|=e_s\le C\left(\sqrt{\f Rm}+\f Rm\right).
\]
This proves \eqref{2-7a}. 

It remains to prove the trace estimate in \eqref{2-6a}. By \eqref{2-6},
\begin{align*}
 |\operatorname{tr}(A_s-I_N)|
 &\le\sum_{i=1}^s|\operatorname{tr}(A_i-A_{i-1})|
            +|\operatorname{tr}(A_0-I_N)|\\
 &\le\f Rm\sum_{i=1}^s2^{-(s-i+1)}+N\|A_0-I_N\|\le\f Rm+\f N{2m}\le\f{3R}{2m}.
\end{align*}

Finally, taking $A=A_s$ and enumerating the points indexed by $J_s$
as $\bx_1,\ldots,\bx_m$, we  complete the proof of the lemma.
\end{proof}
\par\endgroup

\begingroup

\par\endgroup

\subsection{Proof of Theorem \ref{thm-1-1}}

\begingroup

Fix $0<\epsilon<1$ and an orthonormal basis
$\{\phi_j\}_{j=1}^N$ of $X_N\subset L_2(\mu)$, and write
\par\endgroup
\[
 \Phi(\bx):=(\phi_1(\bx),\ldots,\phi_N(\bx))^{\mathsf T},\  \ \bx\in \Og.
\]
\begingroup

Orthonormality and the hypothesis of the theorem give
\[
 \int_\Omega\Phi(\bx)\Phi(\bx)^{\mathsf T}\,d\mu(\bx)=I_N\  \ \text{and}\  \  |\Phi(\bx)|^2\le KN.
\]

\par\endgroup

\begingroup

Every orthogonal transformation of $\Phi$ gives an orthonormal
basis of $X_N$. We shall choose a basis of the form
\par\endgroup
\[
 (\vi_1(\bx),\ldots,\vi_N(\bx))^{\mathsf T}
 =Q^{\mathsf T}\Phi(\bx),
 \qquad Q\in\operatorname{SO}(N).
\]
Thus, for the proof of  {Theorem~\ref{thm-1-1}}, it is enough to find a rotation
$Q\in\operatorname{SO}(N)$ such that, for every integer
$1\le v\le N$, there are
\[
 {m\le CK\epsilon^{-2}v\log(eN/v)}
\]
points $\xi_v^1,\ldots,\xi_v^m\in\Omega$ satisfying
\[
 \sup_{\ta\in\Sigma_v^N\cap\SS^{N-1}}
 \left|
 \ta^{\mathsf T}
 \left(
  \frac 1m\sum_{j=1}^m
       Q^\T\Phi(\xi_v^j)\Phi(\xi_v^j)^{\mathsf T}Q-I_N
 \right)\ta
 \right|
 \le\epsilon.
\]
Here $\Sigma_v^N$ denotes the set of vectors in $\RR^N$
with at most $v$ nonzero coordinates.

For each integer $1\le v\le N$, set
\[
 t_v:=v\log(eN/v)\quad\text{and}\quad
 m_v:=\lceil DK\epsilon^{-2}t_v\rceil,
\]
where $D\ge1$ is an absolute constant to be chosen later.
Throughout the proof, $\|\cdot\|$ denotes the operator norm
induced by the Euclidean norm.

We divide the proof into three steps.

\medskip
\noindent{\bf Step 1. Construction of the matrices.}\\
\begingroup

Applying Lemma~\ref{lem-sampling-matrices} with $m=m_v$ for each
$1\le v\le N$, we obtain positive semidefinite
matrices
\begin{equation}\label{2-8a}
 A_v:=\frac1{m_v}\sum_{j=1}^{m_v}
       \Phi(\by_{v,j})\Phi(\by_{v,j})^{\mathsf T},
 \qquad \by_{v,j}\in\Omega,
 \qquad 1\le v\le N,
\end{equation}
such that
\begin{equation}\label{UDRtrace}
 \left|\frac{\operatorname{tr}A_v}{N}-1\right|
 \le\frac{C_0K}{m_v},
\end{equation}
and 
\begin{equation}\label{2-9b}
 \|A_v-I_N\|\le  \f {C_0}4\Bl(\sqrt{\frac{KN}{m_v}}+\f {KN}{m_v}\Br),
\end{equation}
where $C_0\ge 8$ is an absolute constant. 
In particular, \eqref{2-9b} implies that if $m_v\le KN$, then
\begin{equation}\label{UDRSbounds}
 \|A_v\|\le\frac{C_0KN}{m_v};
\end{equation}
and  if $m_v>KN$, then
\begin{equation}\label{2-9}
 \|A_v-I_N\|\le \delta_v:=C_0\sqrt{\frac{KN}{m_v}}.
\end{equation}
\par\endgroup

\medskip
\noindent{\bf Step 2. Estimates for a random rotation.}\\

Let $Q$ be a random rotation distributed according to normalized Haar measure
on $\operatorname{SO}(N)$. By Lemma~\ref{HaarSphere}, for each
fixed $\ta\in\SS^{N-1}$, the random vector $Q\ta$ is uniformly
distributed on $\SS^{N-1}$. In this step, we prove that 
\beq
 \PP_Q\Bl\{\max_{1\leq v\leq N}\sup_{\ta\in\Sigma_v^N\cap\SS^{N-1}}
 \left|\ta^{\mathsf T}(Q^{\mathsf T}A_vQ-I_N)\ta\right|
 \le\epsilon\Br\}\ge \f 34.\label{UDRallvEvent}
\enq
We first prove that, for each $1\le v\le N$ and every fixed
$\ta\in\SS^{N-1}$,
\beq\label{2-15a}
 \mathbb P_Q\left\{
 \left|\ta^{\mathsf T}Q^{\mathsf T}A_vQ\ta-1\right|
 >\frac{\epsilon}{2}
 \right\}
 \le2\exp(-cDt_v),
\enq
where $c>0$ is an absolute constant.

Suppose first that ${m_v\le KN}$, and define
\[
 F_v(\xi):=\sqrt{\xi^{\mathsf T}A_v\xi}
           =|A_v^{1/2}\xi|,
 \qquad \xi\in\SS^{N-1}.
\]
By \eqref{UDRSbounds}, the Lipschitz constant of $F_v$
is at most
\[
 {L_v:=\sqrt{\frac{C_0KN}{m_v}}.}
\]
Write
\[
 \mathbb E_\sigma h
 :=\int_{\SS^{N-1}}h(\xi)\,d\sigma_N(\xi)\   \ \text{for}\  \  h\in L_1(\SS^{N-1}, \sa_N).
\]
\begingroup

Integrating the concentration estimate in
Lemma~\ref{concentration}, we obtain
\[
 \EE_\sa F_v^2 -(\EE_\sa F_v)^2=\mathbb E_\sigma|F_v-\mathbb E_\sigma F_v|^2
 \le\frac{CL_v^2}{N}
 \le\frac{C_1K}{m_v}.
\]
On the other hand,
\[
 \mathbb E_\sigma F_v^2
 =\int_{\SS^{N-1}}\xi^{\mathsf T}A_v\xi\,d\sigma_N(\xi)
 =\frac{\operatorname{tr}A_v}{N}.
\]
Since $\mathbb E_\sigma F_v\ge0$, it follows from
\eqref{UDRtrace} that
\begin{align*}
 |\mathbb E_\sigma F_v-1|
 &\le |(\mathbb E_\sigma F_v)^2-1|\\
 &\le\left|\frac{\operatorname{tr}A_v}{N}-1\right|
       +\mathbb E_\sigma F_v^2-(\EE_\sa F_v)^2 
 \le{\frac{(C_1+C_0)K}{m_v}}.
\end{align*}
As $m_v\ge DK\epsilon^{-2}t_v$ and $t_v\ge1$,
choosing ${D\ge12(C_1+C_0)}$ gives
\[
 |\mathbb E_\sigma F_v-1|\le\frac{\epsilon}{12}.
\]
\par\endgroup

If $|F_v-1|\le\epsilon/6$, then
\[
 |F_v^2-1|
 \le\frac{\epsilon}{6}
       \left(2+\frac{\epsilon}{6}\right)
 <\frac{\epsilon}{2}.
\]
It follows that
\[
 \left\{|F_v^2-1|>\frac{\epsilon}{2}\right\}
 \subset
 \left\{|F_v-1|>\frac{\epsilon}{6}\right\}
 \subset
 \left\{|F_v-\mathbb E_\sigma F_v|
                  >\frac{\epsilon}{12}\right\}.
\]
Thus,  Lemma~\ref{concentration} implies
\begin{align*}
 \sigma_N\left\{|F_v^2-1|>\frac{\epsilon}{2}\right\}&\leq \sa_N \left\{|F_v-\mathbb E_\sigma F_v|
                  >\frac{\epsilon}{12}\right\}
 \le2\exp\left(-c\frac{N\epsilon^2}{L_v^2}\right)\\
 &\le{2\exp(-c'\epsilon^2m_v/K)}
 \le2\exp(-c_1Dt_v),
\end{align*}
where $c_1>0$ is an absolute constant. Since $Q\ta$
is uniformly distributed on $\SS^{N-1}$, we conclude that
\[
 \mathbb P_Q\left\{
 \left|\ta^{\mathsf T}Q^{\mathsf T}A_vQ\ta-1\right|
 >\frac{\epsilon}{2}
 \right\}
 \le2\exp(-c_1Dt_v).
\]

Now suppose that ${m_v>KN}$, and define
\[
 H_v(\xi):=\xi^{\mathsf T}(A_v-I_N)\xi,
 \qquad \xi\in\SS^{N-1}.
\]
\begingroup

By \eqref{UDRtrace} and the choice of $m_v$, we have
\[
 |\mathbb E_\sigma H_v|
 =\left|\frac{\operatorname{tr}A_v}{N}-1\right|
 \le\frac{C_0K}{m_v}
 \le\frac{C_0\epsilon^2}{Dt_v}
 \le\frac{\epsilon}{4},
\]
provided that $D\ge4C_0$.  Moreover, by \eqref{2-9},
\par\endgroup

\[
 |H_v(\xi)-H_v(\eta)|
 \le2\|A_v-I_N\|\,|\xi-\eta|
 \le2\delta_v|\xi-\eta|,
 \qquad \xi,\eta\in\SS^{N-1}.
\]
\begingroup

Applying Lemma~\ref{concentration} again, we obtain
\begin{align*}
 &\mathbb P_Q\left\{
 \left|\ta^{\mathsf T}Q^{\mathsf T}A_vQ\ta-1\right|
 >\frac{\epsilon}{2}\right\}
 =\sigma_N\left\{|H_v|>\frac{\epsilon}{2}\right\}\\
 &\le\sigma_N\left\{|H_v-\mathbb E_\sigma H_v|
                   >\frac{\epsilon}{4}\right\}
 \le2\exp\left(-c\frac{N\epsilon^2}{\delta_v^2}\right)
 \le2\exp(-c_2Dt_v).
\end{align*}
\par\endgroup
\begingroup

Taking $c=\min\{c_1,c_2\}$ proves \eqref{2-15a}
in both cases.
\par\endgroup

We next pass from a fixed vector to all sparse vectors.
For each $I\subset[N]:=\{1,\ldots,N\}$ with $|I|=v$, set
\[
 {\RR^I:=\{\bx\in\RR^N:\operatorname{supp}\bx\subset I\},
 \qquad \SS_I:=\SS^{N-1}\cap\RR^I.}
\]
Choose a $1/4$-net $\mathcal N_I$ of $\SS_I$ with
$|\mathcal N_I|\le9^v$, and put
\[
 \mathcal N_v
 :=\bigcup_{\substack{I\subset[N]\\|I|=v}}\mathcal N_I.
\]
Then
\[
 |\mathcal N_v|
 \le9^v\binom Nv
 \le\left(\frac{9eN}{v}\right)^v.
\]
The preceding probability estimate \eqref{2-15a} and a union bound give
\begin{align}
 &\mathbb P_Q\left\{
 \max_{\ta\in\mathcal N_v}
 \left|\ta^{\mathsf T}(Q^{\mathsf T}A_vQ-I_N)\ta\right|
 >\frac{\epsilon}{2}
 \right\}\le
 2|\mathcal N_v|\exp(-cDt_v)\notag\\
 &\qquad\le
 2\exp\left(-cDt_v+v\log\frac{9eN}{v}\right)
 \le2e^{-5t_v},\label{UDRnetProbability}
\end{align}
provided that $D$ is sufficiently large. 

To pass from the net to the whole coordinate sphere, fix
$I\subset[N]$ with $|I|=v$. Let $P_I$ be the orthogonal
projection onto $\RR^I$, and define the symmetric operator
\[
 T_I:=P_I (Q^{\mathsf T}A_vQ-I_N) P_I.
\]
Then
\beq\label{UDRcompressionNorm}
 \|T_I\|= \sup_{\bx\in\SS_I}|\bx^{\mathsf T}T_I\bx|
 =\sup_{\bx\in\SS_I}|\bx^{\mathsf T}(Q^{\mathsf T}A_vQ-I_N)\bx|.
\enq
For each $\bx\in\SS_I$, choose $\by\in\mathcal N_I$
with $|\bx-\by|\le1/4$. We have
\begin{equation}\label{2-10}
 |\bx^{\mathsf T}T_I\bx-\by^{\mathsf T}T_I\by|
 \le2\|T_I\|\,|\bx-\by|
 \le\frac12\|T_I\|.
\end{equation}
Taking the supremum over $\bx\in\SS_I$ yields
\[
 \|T_I\|
 \le2\max_{\by\in\mathcal N_I}
          |\by^{\mathsf T}T_I\by|=2 \max_{\by\in\mathcal N_I} |\by^{\mathsf T}(Q^{\mathsf T}A_vQ-I_N)\by| .
\]
\begingroup

By \eqref{UDRcompressionNorm}, it follows that
\begin{align*}
 &\sup_{\ta\in\Sigma_v^N\cap\SS^{N-1}}
 \left|\ta^{\mathsf T}(Q^{\mathsf T}A_vQ-I_N)\ta\right|=\max_{|I|=v}\|T_I\|
 \le2\max_{\ta\in\mathcal N_v}
       \left|\ta^{\mathsf T}(Q^{\mathsf T}A_vQ-I_N)\ta\right|.
\end{align*}
Combining this estimate with \eqref{UDRnetProbability}, we obtain
\begin{align*}
 &\mathbb P_Q\left\{
 \sup_{\ta\in\Sigma_v^N\cap\SS^{N-1}}
 \left|\ta^{\mathsf T}(Q^{\mathsf T}A_vQ-I_N)\ta\right|
 >\epsilon\right\}\\
 &\qquad\le\mathbb P_Q\left\{
 \max_{\ta\in\mathcal N_v}
 \left|\ta^{\mathsf T}(Q^{\mathsf T}A_vQ-I_N)\ta\right|
 >\frac\epsilon2\right\}
 \le2e^{-5t_v}.
\end{align*}
\par\endgroup

\begingroup

Since the function $v\mapsto v\log(eN/v)$ is increasing
on $[1,N]$, we have $t_v\ge\log(eN)$. Another union bound gives
\begin{align*}
 &\mathbb P_Q\left\{
 \max_{1\le v\le N}\sup_{\ta\in\Sigma_v^N\cap\SS^{N-1}}
 \left|\ta^{\mathsf T}(Q^{\mathsf T}A_vQ-I_N)\ta\right|
 >\epsilon\right\}\\
 &\qquad\le\sum_{v=1}^N2e^{-5t_v}
 \le2N(eN)^{-5}<\frac14.
\end{align*}
\par\endgroup
This proves the desired inequality {\eqref{UDRallvEvent}}.\\

\noindent{\bf Step 3. Choice of the orthonormal basis.}\\

\begingroup
By homogeneity  and \eqref{UDRallvEvent}, there exists a rotation
$Q_0\in\operatorname{SO}(N)$ such that
\par\endgroup
\begin{equation}\label{2-12}
 (1-\epsilon)|\ta|^2
 \le\ta^{\mathsf T}Q_0^{\mathsf T}A_vQ_0\ta
 \le(1+\epsilon)|\ta|^2,
 \qquad \ta\in\Sigma_v^N,\  \ 1\leq v\leq N.
\end{equation}
\begingroup

Define an orthonormal basis $\CD_N=\{\vi_j\}_{j=1}^N$
of $X_N$ by
\par\endgroup
\[
 (\vi_1(\bx),\ldots,\vi_N(\bx))^{\mathsf T}
 :=Q_0^{\mathsf T}\Phi(\bx),
 \qquad \bx\in\Omega.
\]

Let $f=\sum_{i=1}^N a_i\vi_i$ with
$\mathbf a\in\Sigma_v^N$. Then
\[
 \|f\|_{L_2(\mu)}^2=|\mathbf a|^2.
\]
Moreover, by \eqref{2-8a},
\begin{align*}
 \frac 1{m_v}\sum_{j=1}^{m_v}
       |f(\by_{v,j})|^2
 &=
 \frac 1{m_v}\sum_{j=1}^{m_v}
       |\mathbf a^{\mathsf T}Q_0^{\mathsf T}\Phi(\by_{v,j})|^2{=\mathbf a^{\mathsf T}Q_0^{\mathsf T}A_vQ_0\mathbf a.}
\end{align*}
Applying \eqref{2-12}, we obtain
\[
 (1-\epsilon)\|f\|_{L_2(\mu)}^2
 \le
 \frac 1{m_v}\sum_{j=1}^{m_v}
       |f(\by_{v,j})|^2
 \le
 (1+\epsilon)\|f\|_{L_2(\mu)}^2.
\]
Finally, since $D$ is now a fixed absolute constant,
\begingroup

\[
 m_v=\lceil DK\epsilon^{-2}t_v\rceil
 \le(D+1)K\epsilon^{-2}t_v
 \le CK\epsilon^{-2}v\log(eN/v),
 \qquad 1\le v\le N.
\]
\par\endgroup

Thus, the same basis works for every $1\le v\le N$.
This completes the proof of Theorem~\ref{thm-1-1}.

\subsection{Proof of Corollary \ref{UDRchosen}}

\begingroup

Fix an orthonormal basis $\{\phi_j\}_{j=1}^N$ of $X_N$, and set
\[
 K_N(\bx):=\sum_{j=1}^N|\phi_j(\bx)|^2,
 \qquad w(\bx):=\frac{K_N(\bx)}N
\]
and $\Omega_+=\{\bx\in\Omega:K_N(\bx)>0\}$. Then every
function in $X_N$ vanishes on $\Omega\setminus\Omega_+$.
Furthermore, 
\[
 \lambda_N(\bx)=\frac1{K_N(\bx)},
 \qquad \bx\in\Omega_+.
\]

Define a measure $\nu$ on $\Omega_+$ by $d\nu=w\,d\mu$.
Since the functions $\phi_j$ are orthonormal,
\[
 \int_{\Omega_+}w\,d\mu
 =\frac1N\sum_{j=1}^N\|\phi_j\|_2^2=1,
\]
so $\nu$ is a probability measure. Put
\[
 \psi_j(\bx):=w(\bx)^{-1/2}\phi_j(\bx),
 \qquad \bx\in\Omega_+,
 \qquad 1\le j\le N,
\]
and let
$\widetilde X_N:=\operatorname{span}\{\psi_1,\ldots,\psi_N\}$.
For every $f\in X_N$, we have
\[
 \|w^{-1/2}f\|_{L_2(\nu)}^2
 =\int_{\Omega_+}|f|^2\,d\mu
 =\|f\|_{L_2(\mu)}^2.
\]
Thus $f\mapsto w^{-1/2}f$ is a linear isometry from $X_N$
onto $\widetilde X_N$, and $\{\psi_j\}_{j=1}^N$ is an
orthonormal basis of $\widetilde X_N$ in $L_2(\nu)$. Moreover,
\begin{equation}\label{OBlength}
 \sum_{j=1}^N|\psi_j(\bx)|^2=N,
 \qquad \bx\in\Omega_+.
\end{equation}
Thus,  $\widetilde X_N$ consists of bounded measurable
functions and satisfies the hypothesis of Theorem~\ref{thm-1-1}
with $K=1$.

Applying Theorem~\ref{thm-1-1} to
$\widetilde X_N$ on $(\Omega_+,\nu)$, we obtain an orthonormal
basis $\{\widetilde\vi_j\}_{j=1}^N$, independent of $v$, such
that, for every $1\le v\le N$, there are points
$\xi^1,\ldots,\xi^m\in\Omega_+$ with
\[
 m\le C\epsilon^{-2}v\log(eN/v)
\]
and
\[
 (1-\epsilon)\|g\|_{L_2(\nu)}^2
 \le\frac1m\sum_{i=1}^m|g(\xi^i)|^2
 \le(1+\epsilon)\|g\|_{L_2(\nu)}^2,
 \qquad
 g\in\Sigma_v(\{\widetilde\vi_j\}_{j=1}^N).
\]
Define
\[
 \vi_j(\bx):=
 \begin{cases}
  w(\bx)^{1/2}\widetilde\vi_j(\bx),&\bx\in\Omega_+,\\
  0,&\bx\in\Omega\setminus\Omega_+.
 \end{cases}
\]
By the isometry above, $\CD_N:=\{\vi_j\}_{j=1}^N$ is an
orthonormal basis of $X_N$. If
$f=\sum_{j\in I}a_j\vi_j$ with $|I|\le v$, then
\[
 g:=w^{-1/2}f=\sum_{j\in I}a_j\widetilde\vi_j
 \quad\text{on }\Omega_+,
 \qquad
 \|g\|_{L_2(\nu)}=\|f\|_{L_2(\mu)}.
\]
Since $w^{-1}=N\lambda_N$ on $\Omega_+$, we also have
\[
 \frac1m\sum_{i=1}^m|g(\xi^i)|^2
 =\frac Nm\sum_{i=1}^m\lambda_N(\xi^i)|f(\xi^i)|^2.
\]
{Substituting these identities into the discretization
inequality for $g$, we obtain}
\[
 (1-\epsilon)\|f\|_{L_2(\mu)}^2
 \le\frac Nm\sum_{i=1}^m\lambda_N(\xi^i)|f(\xi^i)|^2
 \le(1+\epsilon)\|f\|_{L_2(\mu)}^2,
 \qquad f\in\Sigma_v(\CD_N).
\]

\par\endgroup


\section{Random orthonormal bases}\label{OB}
Our main purpose in this section is to prove simultaneous weighted
$L_p$ estimates for a random orthonormal basis and a logarithmic
bound for the integral of its squared pointwise maximum.
Let $(\Omega,\mu)$ be a probability space, and let $X_N$ be a real
$N$-dimensional subspace of $L_2(\mu)$, $N\ge2$. We use the
measurable representatives fixed in the introduction, so every
function in $X_N$ is defined everywhere on $\Omega$.
Let $\lambda_N$ be the Christoffel function associated with
$X_N$ and $\mu$, as defined in \eqref{2-1}.
Fix an orthonormal basis $\{\phi_j\}_{j=1}^N$ of
$(X_N,\|\cdot\|_{L_2(\mu)})$, and put
\[
 \Phi(\bx):=(\phi_1(\bx),\ldots,\phi_N(\bx))^{\mathsf T},
 \qquad \bx\in\Omega.
\]
For each rotation $Q\in\operatorname{SO}(N)$, define
\[
 \Phi_Q(\bx):=Q\Phi(\bx)
 =\bigl(\phi_{Q,1}(\bx),\ldots,\phi_{Q,N}(\bx)\bigr)^{\mathsf T},
 \qquad \bx\in\Omega.
\]
Then $\{\phi_{Q,j}\}_{j=1}^N$ is also an orthonormal basis
of $(X_N,\|\cdot\|_{L_2(\mu)})$.
Set $\Omega_+:=\{\bx\in\Omega:\Phi(\bx)\ne0\}$.
Recall that $\lambda_N=|\Phi|^{-2}$ on $\Omega_+$ and
$\lambda_N=\infty$ on $\Omega\setminus\Omega_+$.
We use this notation throughout the section.

\begin{Theorem}\label{OBT1}
Let $Q$ be an $N\times N$ random matrix distributed according
to normalized Haar measure on $\operatorname{SO}(N)$, and let
$\{\phi_{Q,j}\}_{j=1}^N$ be the random orthonormal basis of
$X_N$ defined above. Then there is an absolute constant $C>0$
such that, for every $\tau\ge1$, the inequalities
\begin{equation}\label{OB0}
 \max_{1\le j\le N}
 \left(\int_\Omega |\phi_{Q,j}(\bx)|^p
       (N\lambda_N(\bx))^{\f p2-1}\,d\mu(\bx)\right)^{\f1p}
 \le C\tau\sqrt p,
 \qquad 1\le p<\infty,
\end{equation}
hold simultaneously with probability at least
$1-2e^{-\tau^2}N^{1-8\tau^2}$.
Here and below, the weighted integrands are defined to be zero on $\Omega\setminus\Omega_+$.
\end{Theorem}

\begin{Corollary}\label{UDRallvLq}
For each fixed $0<\epsilon<1$, the orthonormal basis
$\CD_N=\{\vi_j\}_{j=1}^N$ in Corollary~\ref{UDRchosen} can be
chosen so that, in addition to the discretization inequalities
for every $1\le v\le N$,
\beq\label{3-9}
 \max_{1\le j\le N}
 \left(\int_\Omega |\vi_j(\bx)|^p
       (N\lambda_N(\bx))^{\frac{p}{2}-1}\,d\mu(\bx)\right)^{\frac{1}{p}}
 \le C\sqrt p,
 \qquad 1\le p<\infty.
\enq
Here $C$ is an absolute constant.
\end{Corollary}
The following corollary gives a logarithmic bound for the $L_2(\mu)$ norm of the pointwise maximum of a basis satisfying \eqref{OB0}.

\begin{Corollary}\label{UDRlogenv}
Let $Q\in\operatorname{SO}(N)$ be a rotation for which \eqref{OB0} holds for some $\tau\ge1$ and every $1\le p<\infty$. Then there is an absolute constant $C>0$ such that
\begin{equation}\label{UDRenvelope}
 \left(\int_\Omega\max_{1\le j\le N}
                       |\phi_{Q,j}(\bx)|^2\,d\mu(\bx)\right)^{\frac{1}{2}}\le C\tau\sqrt{\log(eN)}.
\end{equation}
\end{Corollary}

The logarithmic order in \eqref{UDRenvelope} cannot be improved
for arbitrary $N$-dimensional subspaces.
\begin{Proposition}\label{UDRsharpEnvelope}
For every $N\ge2$, there are a probability space $(\Omega,\mu)$ and a real $N$-dimensional subspace $X_N$ of $L_2(\mu)$ such that every orthonormal basis $\{\vi_j\}_{j=1}^N$ of $X_N$ satisfies
\[
 \Bl(\int_{\Omega}\max_{1\le j\le N}|\vi_j(\bx)|^2\,d\mu(\bx)\Br)^{\frac{1}{2}}\ge c\sqrt{\log(eN)},
\]
where $c>0$ is an absolute constant.
\end{Proposition}
\begin{Remark}
{
The change-of-density argument below follows the sampling
construction studied by Adcock, Cardenas, Dexter, and Moraga
\cite[Section~4.4]{ACDM2022}. For a prescribed finite dictionary,
these authors choose a sampling density proportional to the
pointwise maximum of the squared dictionary elements. They show
that this choice minimizes the weighted uniform bound in their
sampling estimates. Together with their recovery estimates for
$\ell_1$-minimization, this gives stable sparse approximation from
independent samples. In the orthonormal case, the resulting sample
bound is proportional to the sparsity times the integral of this
pointwise maximum, up to logarithmic factors. Thus, their construction
applies directly to the bases obtained in Corollary~\ref{UDRlogenv}.
For a detailed account of the related Christoffel sampling method
for weighted least-squares approximation, with sample complexity
of order $N\log N$, see Adcock \cite[Sections~5 and~6]{Adcock2025}.

To make this connection explicit, let $Q\in\operatorname{SO}(N)$
be a rotation for which
\beq\label{3-3b}
 A_N:=\int_\Omega\max_{1\le j\le N}
                 |\phi_{Q,j}(\bx)|^2\,d\mu(\bx)
 \le C\log(eN).
\enq
Define
\begin{equation}\label{UDRsecondmeasure}
 \rho(\bx):=
 \frac{1+\max_{1\le j\le N}|\phi_{Q,j}(\bx)|^2}{1+A_N},
 \qquad \bx\in\Omega.
\end{equation}
Clearly,  $\rho$ is strictly positive on  $\Omega$,
and $\int_\Omega\rho\,d\mu=1$. Put
\[
 \vi_j(\bx):=\frac{\phi_{Q,j}(\bx)}{\sqrt{\rho(\bx)}},
 \qquad \bx\in\Omega,\quad 1\le j\le N.
\]
Then $\{\vi_j\}_{j=1}^N$ is orthonormal in $L_2(\rho\,d\mu)$, and satisfies 
\[
 |\vi_j(\bx)|^2
 =\frac{|\phi_{Q,j}(\bx)|^2}{\rho(\bx)}\le1+A_N,
 \qquad \bx\in\Omega.
\]
On the other hand,
\[
 A_N=\int_\Omega\max_{1\le j\le N}|\vi_j(\bx)|^2
                         \rho(\bx)\,d\mu(\bx)
 \le\max_{1\le j\le N}\|\vi_j\|_{L_\infty(\rho\,d\mu)}^2.
\]
Consequently,
\[
 A_N\le\max_{1\le j\le N}\|\vi_j\|_{L_\infty(\rho\,d\mu)}^2
 \le1+A_N.
\]
In particular, \eqref{3-3b} and the pointwise bound above give
\begin{equation}\label{UDRbounded}
 \max_{1\le j\le N}\sup_{\bx\in\Omega}|\vi_j(\bx)|^2
 \le C\log(eN).
\end{equation}
Thus Corollary~\ref{UDRlogenv} supplies a choice of orthonormal
basis for which the sampling parameter in \cite{ACDM2022}
is  bounded by \(C\log(eN)\).
\par}
\end{Remark}

\subsection{Proof of Theorem \ref{OBT1}}\label{subsection-3-1}
We first use a change of density to reduce the weighted
estimate in \eqref{OB0} to an $L_p$ estimate for a normalized
system. Put
\[
 K_N(\bx):=|\Phi(\bx)|^2,\qquad w(\bx):=\frac{K_N(\bx)}N=(N\ld_N(\bx))^{-1}.
\]
Define $ d\nu:=w\,d\mu$.
Since $\int_\Omega K_N\,d\mu=N$, the measure $\nu$ is a
probability measure with $\nu(\Omega_+)=1$. Define
$\psi_j=w^{-1/2}\phi_j$ on $\Omega_+$ and $\psi_j=0$
elsewhere. Then $\{\psi_j\}_{j=1}^N$ is an orthonormal basis of
\beq 
 \widetilde X_N:=\operatorname{span}\{\psi_1,\ldots,\psi_N\}
 \subset L_2(\nu),\label{3-3-a}
\enq
and the vector $\Psi:=(\psi_1,\ldots,\psi_N)^{\mathsf T}$
satisfies $|\Psi|^2=N$ on $\Omega_+$.
For each rotation $Q$, write
\beq\label{3-4-a}
 \Psi_Q:=Q\Psi=(\psi_{Q,1},\ldots,\psi_{Q,N})^{\mathsf T}.
\enq
Since $\psi_{Q,j}=w^{-1/2}\phi_{Q,j}$ and $w^{-1}=N\lambda_N$ on $\Omega_+$, we have
\begin{equation}\label{OBweightedidentity}
 \|\psi_{Q,j}\|_{L_p(\nu)}^p
 =\int_{\Omega_+}|\phi_{Q,j}|^p w^{1-\f p2}\,d\mu
 =\int_\Omega|\phi_{Q,j}|^p(N\lambda_N)^{\f p2-1}\,d\mu,
 \qquad 1\le p<\infty.
\end{equation}
Thus it suffices to estimate the $L_p(\nu)$ norms of the
functions $\psi_{Q,j}$. We need the following consequence
of spherical concentration.

\begin{Lemma}\label{OBL1}
There exists an absolute constant $C>0$ such that, for every integer
$N\ge2$, every $2\le q<\infty$, and every $\tau\ge1$,
\beq\label{3-2a}
 \sigma_N\left\{\bu\in\SS^{N-1}:
  \left\|\sum_{j=1}^N u_j\psi_j\right\|_{L_q(\nu)}
  >C\tau\sqrt q\right\}
 \le2\exp(-4\tau^2qN^{2/q}).
\enq
\end{Lemma}

\begin{proof}
 By \eqref{OBlength} and the
Cauchy--Schwarz inequality,
\[
 \|f\|_{L_\infty(\nu)}\le\sqrt N\|f\|_{L_2(\nu)},
 \qquad f\in\widetilde X_N.
\]
Consequently, for $q\ge 2$,
\beq\label{3-2}
 \|f\|_{L_q(\nu)}
 \le\|f\|_{L_2(\nu)}^{2/q}\|f\|_{L_\infty(\nu)}^{1-2/q}
 \le N^{\f12-\f1q}\|f\|_{L_2(\nu)},
 \qquad f\in\widetilde X_N.
\enq
Consider the function
\[
 F_q(\bu):=\|\langle\bu,\Psi\rangle\|_{L_q(\nu)},
 \qquad \bu\in\SS^{N-1},
\]
where $\langle\cdot,\cdot\rangle$ denotes the Euclidean inner
product on $\RR^N$. By \eqref{3-2} and orthonormality,
\begin{align*}
 |F_q(\bu)-F_q(\bv)|
 &\le\left\|\sum_{j=1}^N(u_j-v_j)\psi_j\right\|_{L_q(\nu)}\\
 &\le N^{\f 12-\f1q}
       \left\|\sum_{j=1}^N(u_j-v_j)\psi_j\right\|_{L_2(\nu)}
 =N^{\f 12-\f1q}|\bu-\bv|.
\end{align*}
Thus $F_q:\SS^{N-1}\to\RR$ has Lipschitz constant at most $N^{\f12-\f1q}$. Applying Lemma~\ref{concentration},
we obtain
\begin{equation}\label{OB8}
 \sigma_N\Bl\{|F_q-\mathbb E_\sigma F_q|>C_1N^{-1/q}t\Br\}
 \le2\exp(-4t^2),\qquad t>0,
\end{equation}
where $C_1$ is an absolute constant and $\mathbb E_\sigma$
denotes integration with respect to the measure  $\sigma_N$ on $\SS^{N-1}$.

By H\"older's inequality, Fubini's theorem, and rotational
invariance of $\sigma_N$,
\begin{align*}
 \mathbb E_\sigma F_q
 &\le(\mathbb E_\sigma F_q^q)^{1/q}=\left(\int_\Omega\int_{\SS^{N-1}}
       |\langle\bu,\Psi(\bx)\rangle|^q
       \,d\sigma_N(\bu)\,d\nu(\bx)\right)^{1/q}\\
 &=\sqrt N\left(\int_{\SS^{N-1}}|u_1|^q\,d\sigma_N(\bu)\right)^{1/q}
 \le C_2\sqrt q.
\end{align*}
Taking $t=\tau\sqrt q\,N^{1/q}$ in \eqref{OB8} and
$C=C_1+C_2$, we obtain \eqref{3-2a}, since $\tau\ge1$.
\end{proof}

We now complete the proof of Theorem~\ref{OBT1}.

For every rotation $Q\in\operatorname{SO}(N)$, the system $\{\psi_{Q,j}\}_{j=1}^N$
is orthonormal in $L_2(\nu)$.
For each $1\le j\le N$, let $\mathbf e_j$ be the $j$th
coordinate vector. Then
\[
 \psi_{Q,j}=\langle Q^{\mathsf T}\mathbf e_j,\Psi\rangle.
\]
By Lemma~\ref{HaarSphere} and the invariance of Haar measure
under transposition, $Q^{\mathsf T}\mathbf e_j$ is uniformly
distributed on $\SS^{N-1}$. Lemma~\ref{OBL1} and a union
bound therefore give, for each integer $k\ge3$,
\eq{
 \mathbb P\left\{\max_{1\le j\le N}
       \|\psi_{Q,j}\|_{L_k(\nu)}>C_0\tau\sqrt k\right\}      &\le N\sigma_N\left\{\bu\in\SS^{N-1}:\|\langle\bu,\Psi\rangle\|_{L_k(\nu)}>C_0\tau\sqrt k\right\}\\
 &\le2N\exp(-4\tau^2kN^{2/k}).
}
By  monotonicity of the $L_p(\nu)$ norms, we can find an absolute constant $C>4C_0$ such that
\begin{align*}
 \mathbb P\left\{
  \sup_{1\leq p<\infty}\max_{1\le j\le N}
  \frac{\|\psi_{Q,j}\|_{L_p(\nu)}}{\sqrt p}>C\tau\right\}&\leq 
  \mathbb P\left\{
  \sup_{k\in \mathbb{N}, k\ge 3}\max_{1\le j\le N}
  \frac{\|\psi_{Q,j}\|_{L_k(\nu)}}{\sqrt k}>C_0\tau\right\}\\
 &\le2N\sum_{k=3}^\infty e^{-4\tau^2kN^{2/k}}.
\end{align*}
Using $e^x\ge x+1$, we have 
\[
 kN^{2/k}=k\exp\Bl(\f{2\log N}k\Br)\ge k+2\log N, \  \ \forall k\ge 3,
\]
and hence 
\begin{align*}
 2N\sum_{k=3}^\infty e^{-4\tau^2kN^{2/k}}
 &\le2N^{1-8\tau^2}\sum_{k=3}^\infty e^{-4\tau^2k}
 \le2e^{-\tau^2}N^{1-8\tau^2}.
\end{align*}
 Combining  the preceding estimate with \eqref{OBweightedidentity} proves the theorem.

\subsection{Proof of Corollaries \ref{UDRallvLq} and   \ref{UDRlogenv}}

\begin{proof}[Proof of Corollary \ref{UDRallvLq}] Let $\wt X_N$ denote the $N$-dimensional space defined by \eqref{3-3-a}.
We apply the construction in the proof of
Theorem~\ref{thm-1-1} to $\widetilde X_N$ with $K=1$.
By \eqref{UDRallvEvent} and the invariance of Haar measure under transposition, there is an event $\mathcal D_\epsilon$
of probability at least $3/4$ on which the basis $Q\Psi$ given in \eqref{3-4-a}
satisfies the required discretization inequalities in Theorem~\ref{thm-1-1},
simultaneously for every $v$.
The change of density in the proof of Corollary~\ref{UDRchosen}
transfers these inequalities to the basis $Q\Phi$ in $L_2(\mu)$.
Let $\mathcal E_1$ be the event on which \eqref{OB0} holds
with $\tau=1$ for every $1\le p<\infty$.
By Theorem~\ref{OBT1},
\[
 \mathbb P(\mathcal D_\epsilon\cap\mathcal E_1)
 \ge\frac34-2e^{-1}N^{-7}>0.
\]
Choose $Q$ in this intersection and set $\vi_j=\phi_{Q,j}$.
Then \eqref{UDRdiscretization} and the asserted weighted
$L_p(\mu)$ estimates \eqref{3-9} hold for the same basis.
\end{proof}
\begin{proof}[Proof of Corollary  \ref{UDRlogenv}]
 We use the notation of Section~\ref{subsection-3-1}. Fix $\tau\ge1$, and let $Q\in\operatorname{SO}(N)$ be a rotation for which
\eqref{OB0} holds for every $1\le p<\infty$; that is,
\[ \Bl(\max_{1\leq j\leq N} \int_{\Og} |\phi_{Q,j} (\bx)|^p w(\bx)^{1-\f p2} \, d\mu(\bx) \Br)^{\f 1p}\leq C \tau \sqrt{p},\  \ \forall 1\leq p<\infty,\]
where  $w(\bx)=(N\ld_N(\bx))^{-1}$.  
Put
\[
 A_N(Q):=\int_\Omega\max_{1\le j\le N}
                       |\phi_{Q,j}(\bx)|^2\,d\mu(\bx).\]
Applying H\"older's inequality, for every $p>2$ we obtain
\eq{A_N(Q)&\leq \int_{\Og} \Bl( \sum_{j=1}^N |\phi_{Q,j}(\bx)|^pw(\bx)^{1-\f p2}\Br)^{\f 2p}w(\bx)^{1-\f 2p}\, d\mu(\bx)\\
&\leq \left(\sum_{j=1}^N \int_{\Og}  |\phi_{Q,j}(\bx)|^pw(\bx)^{1-\f p2}\, d\mu(\bx)\right)^{\f 2p} \Bl(\int_{\Og} w(\bx)\, d\mu(\bx)\Br)^{1-\f2p}\\
&\leq C \tau^2 p N^{\f 2p}.
}
Taking $p=3\log (4N)$, we deduce
\[
 A_N(Q)
 \le C\tau^2\log(eN).
\]
This proves the corollary.
\end{proof}

\subsection{Proof of Proposition \ref{UDRsharpEnvelope}}

Let $(\Omega,\mu)=(\SS^{N-1},\sigma_N)$. For $1\le j\le N$, define $\phi_j:\SS^{N-1}\to\RR$ by

\[ \phi_j(\bx)= \sqrt N x_j,\   \ \bx\in\SS^{N-1}.\]
Then $\{\phi_j\}_{j=1}^N$ is an orthonormal system in $L_2(\SS^{N-1},\sigma_N)$, and 
$$\Phi(\bx):=(\phi_1(\bx), \ldots, \phi_N(\bx))^\T=\sqrt N\,\bx,\  \ \bx\in \SS^{N-1}.$$
Let
\[
 X_N:=\operatorname{span}\{\phi_j:\  1\le j\le N\}.
\]
Let $\{u_j\}_{j=1}^N$ be an arbitrary orthonormal basis of $X_N$. Then there is an orthogonal matrix $Q\in\operatorname{O}(N)$ such that

\[ U(\bx):=(u_1(\bx), \ldots, u_N(\bx))^\T=Q \Phi(\bx)=\sqrt{N} Q \bx,\  \ \bx\in \SS^{N-1}.\]
 By invariance of $\sigma_N$ under orthogonal transformations, we obtain

\eq{
 \int_{\SS^{N-1}}&\max_{1\leq j\leq N}|u_j(\bx)|^2\,d\sigma_N(\bx)
   =N\int_{\SS^{N-1}}\max_{1\leq j\leq N} |\la Q \bx, \mathbf{e}_j\ra|^2\,d\sigma_N(\bx)\\
 &=N\int_{\SS^{N-1}} \max_{1\leq j\leq N} |x_j|^2\, d\sa_N(\bx)
 \asymp \log (eN).
}
This proves the proposition. 



\section{Changes of measure and universal discretization}
\label{Ex}\label{UD}
We now consider sampling in a prescribed system of functions. An
integrable pointwise bound permits a change of measure that makes the system
uniformly bounded and preserves its coefficients. We first show that this change of measure preserves
recovery errors in $L_p$. We then give the quadratic
discretization estimate used in Section~\ref{SR}.

\subsection{A change of measure in \texorpdfstring{$L_p$}{Lp}}
Let $1\le p<\infty$, and let $\Phi=\{\phi_\bk\}_{\bk\in\Lambda}$ be a
finite or countable system defined everywhere on $\Og$, where $\Lambda\subset\bbZ^d$. Suppose that the chosen pointwise representatives satisfy
\be\label{Lp1}
 |\phi_\bk(\bx)|\le B(\bx),\   \   \bx\in\Omega,\quad \bk\in\Lambda,
\ee
for some function  $B\in L_p(\mu)$ satisfying $1\le B(\bx)<\infty$ for all $\bx\in\Og$.
 Define $c_p>0$ and $\alpha_p$ by
\be\label{Lp2}
 c_p^p\int_\Omega B(\bx)^p\,d\mu=1,
 \qquad \alpha_p(\bx):=(c_pB(\bx))^{-1},
\ee
and put
\be\label{Lp3}
 U_pf:=\alpha_pf,\qquad \psi_\bk:=U_p\phi_\bk,\qquad
 d\nu_p:=c_p^pB^p\,d\mu.
\ee
 The measure
$\nu_p$ is a probability measure on $\Og$, and
\be\label{Lp8}
 |\psi_\bk(\bx)|\le c_p^{-1},\qquad \bx\in\Omega.
\ee
For every $f\in L_p(\mu)$, the definition of the density gives
\be\label{Lp9}
 \|U_pf\|_{L_p(\nu_p)}^p
 =\int_\Omega |f|^p\alpha_p^pc_p^pB^p\,d\mu
 =\|f\|_{L_p(\mu)}^p.
\ee

For a finite dictionary $\cD_\Gamma=\{g_\bk\}_{\bk\in \Gamma}$, write
$[\cD_\Gamma]=\operatorname{span}\cD_\Gamma$. We denote by
$\varrho_m^o(\bF,\cD_\Gamma,L_p(\eta))$ the recovery error defined with respect to a measure $\eta$ in
Section~\ref{In}, with reconstruction maps restricted to take values
in $[\cD_\Gamma]$. Put $U_p\bF:=\{U_pf:f\in\bF\}$, and write
$\Phi_\Gamma=\{\phi_\bk\}_{\bk\in \Gamma}$ and
$\Psi_\Gamma=\{\psi_\bk\}_{\bk\in \Gamma}$ for finite $\Gamma\subset\Lambda$.

\begin{Proposition}\label{LpRecoveryTransport}\label{SRL1}
Under the preceding assumptions, for every class
$\bF\subset L_p(\mu)$ with specified pointwise representatives and
every finite $\Gamma\subset\Lambda$,
\[
 \varrho_m^o(\bF,\Phi_\Gamma,L_p(\mu))
 =\varrho_m^o(U_p\bF,\Psi_\Gamma,L_p(\nu_p)).
\]
The same equality holds when the output space is unrestricted.
\end{Proposition}
\begin{proof}
Fix sampling points $\xi^1,\ldots,\xi^m\in\Omega$, and let $D_\xi$
be the invertible diagonal matrix with entries $\alpha_p(\xi^i)$. The sample
vector of $U_pf$ is $D_\xi$ times the sample vector of $f$. A
reconstruction map $\cM$ for the transformed class induces
$U_p^{-1}\cM D_\xi$ for the original class, with the same error by
\eqref{Lp9}. Conversely, a reconstruction map for the original class
induces $U_p\cM D_\xi^{-1}$ for the transformed class. Since $U_p$ maps
$[\Phi_\Gamma]$ onto $[\Psi_\Gamma]$, this correspondence preserves the output
restrictions. Taking the infimum over reconstruction maps and sampling
points proves the result. Unique dictionary representations are not
required.
\end{proof}

\subsection{The \texorpdfstring{$L_2$}{L2} specialization}
Suppose now that $\Phi=\{\phi_\bk\}_{\bk\in\bbZ^d}\subset L_2(\mu)$
has a common measurable pointwise bound $B\in L_2(\mu)$ satisfying
\be\label{In1}
 |\phi_\bk(\bx)|\le B(\bx),\qquad
 1\le B(\bx)<\infty,\quad \bx\in\Omega,\quad \bk\in\bbZ^d.
\ee
We use the normalization
\be\label{In2}
 c_0:=\|B\|_2^{-1},\qquad c_0^2\int_\Omega B^2\,d\mu=1.
\ee
This assumption includes systems which are not uniformly bounded. For
example, the orthonormal Legendre polynomials on $[-1,1]$, with normalized
Lebesgue measure, are bounded pointwise by $C(1-x^2)^{-1/4}$ on $(-1,1)$.

For $p=2$, use the notation $c_0=\|B\|_{L_2(\mu)}^{-1}$ from
\eqref{In2}, and put $\alpha=(c_0B)^{-1}$ and $\nu=\nu_2$. Then
\be\label{In3}
 \psi_\bk=\alpha\phi_\bk,\qquad d\nu=c_0^2B^2\,d\mu.
\ee
The transformation preserves every Gram matrix:
\be\label{In4}
 \langle\psi_\bk,\psi_{\bk'}\rangle_\nu
 =\int_\Omega\phi_\bk\,\overline{\phi_{\bk'}}\,\alpha^2c_0^2B^2\,d\mu
 =\langle\phi_\bk,\phi_{\bk'}\rangle_\mu.
\ee
In particular, it preserves orthonormality. Moreover,
\be\label{In5}
 |\psi_\bk(\bx)|\le c_0^{-1},\qquad \bx\in\Omega,
\ee
and, for every $g\in L_2(\mu)$,
\be\label{In6}
 \|\alpha g\|_{L_2(\nu)}=\|g\|_{L_2(\mu)}.
\ee

\begin{Proposition}\label{InP1}\label{InP2}
The systems $\Phi$ in $L_2(\mu)$ and $\Psi$ in $L_2(\nu)$ have the
same lower Riesz constants. More generally, if
\begin{equation}\label{Riesz}
 R_1\left(\sum_\bk|a_\bk|^2\right)^{1/2}
 \le\left\|\sum_\bk a_\bk\phi_\bk\right\|_{L_2(\mu)}
 \le R_2\left(\sum_\bk|a_\bk|^2\right)^{1/2}
\end{equation}
for all finite coefficient vectors, where $0<R_1\le R_2<\infty$, then
the same inequalities hold for $\Psi$ in $L_2(\nu)$ with the same
constants $R_1,R_2$.
\end{Proposition}
\begin{proof}
For every finitely supported coefficient sequence, \eqref{In6} gives
\[
 \left\|\sum_\bk a_\bk\phi_\bk\right\|_{L_2(\mu)}
 =\left\|\sum_\bk a_\bk\psi_\bk\right\|_{L_2(\nu)}.
\]
The assertion follows by applying this identity to the
Riesz inequalities.
\end{proof}

\subsection{Quadratic discretization}
\begin{Definition}\label{UDD1}
Let $1\le p<\infty$, and let $\cX=\{X(n)\}_{n=1}^k$ be a collection
of finite-dimensional subspaces of $L_p(\mu)$ with specified
representatives. Points $\xi^1,\ldots,\xi^m\in\Omega$ and
nonnegative weights $w_1,\ldots,w_m$ provide \emph{$L_p(\mu)$-universal
weighted sampling discretization} with constants $0<C_1\le C_2$ if
\be\label{ud}
 C_1\|f\|_p^p\le\sum_{i=1}^m w_i|f(\xi^i)|^p
 \le C_2\|f\|_p^p,\qquad f\in\bigcup_{n=1}^k X(n).
\ee
Repetitions of the sampling points are allowed. We denote the least
possible number of points by $m(\cX,p,C_1,C_2)$ and use the abbreviation
$L_p$-uwsd$(C_1,C_2)$ for this property.
\end{Definition}

We apply this definition to the coordinate subspaces whose union is
$\Sigma_s(\cD_N)$. The following consequence of
\cite[Theorem~1.1]{BDJR} requires only a lower Riesz inequality.

\begin{Theorem}\label{UDTsharp}
Let $\eta$ be a probability measure on $\Og$, let
$\cD_N=\{g_j\}_{j=1}^N$ satisfy $|g_j(\bx)|\le L$ on $\Omega$, and let
$1\le s\le N$. Suppose that for every coefficient vector supported on
at most $s$ indices,
\begin{equation}\label{UDsharplower}
 \sum_{j=1}^N|a_j|^2\le K
 \left\|\sum_{j=1}^N a_jg_j\right\|_{L_2(\eta)}^2.
\end{equation}
Put $R=KL^2$. There are absolute constants $C,c>0$ such that, whenever
\begin{equation}\label{UDsharpm}
 m\ge C R s\log(eN)\log^2(eRs),
\end{equation}
independent points $\xi^1,\ldots,\xi^m$ drawn from $\eta$ satisfy
\begin{equation}\label{UDsharpineq}
 \frac12\|f\|_{L_2(\eta)}^2\le
 \frac1m\sum_{i=1}^m|f(\xi^i)|^2
 \le\frac32\|f\|_{L_2(\eta)}^2,
 \qquad f\in\Sigma_s(\cD_N),
\end{equation}
with probability at least $1-2\exp(-cm/(Rs))$.
\end{Theorem}
\begin{proof}
Applying \eqref{UDsharplower} to a single coordinate gives $R\ge1$.
Consider the normalized coefficient set
\[
 T:=\left\{a\in\bbC^N:\bigl|\{j:a_j\ne0\}\bigr|\le s,
 \ \left\|\sum_{j=1}^N a_jg_j\right\|_{L_2(\eta)}=1\right\}.
\]
For $a\in T$, Cauchy--Schwarz and \eqref{UDsharplower} give
\[
\|a\|_1\le\sqrt{s}\|a\|_2\le\sqrt{Ks}.
\]
Let
$X=(\overline{g_j(\xi)}/L)_{j=1}^N$, where $\xi$ has distribution
$\eta$, and use the Hermitian pairing that is linear in its first
argument. The coordinates of $X$ have absolute value at most one, while
\[
 LT\subset\sqrt{Rs}\,B_1^N,\qquad
 \mathbb E|\langle La,X\rangle|^2=1\quad(a\in T),
\]
where $B_1^N=\{u\in\bbC^N:\|u\|_1\le1\}$.

Apply \cite[Theorem~1.1]{BDJR} to $LT$, with coordinate bound one and
size parameter $Rs$ (or its ceiling). In that theorem the deviation is
bounded by $2c_1\delta$, since the quadratic expectation equals one on
$LT$. Fix an absolute $\delta\in(0,\kappa)$ with
$2c_1\delta\le1/2$. The required lower bound on $m$ follows from
\eqref{UDsharpm}, and its failure probability is at most
$2\exp(-cm/(Rs))$ after changing absolute constants. This proves
\eqref{UDsharpineq} on the normalized set $T$. 

Homogeneity gives the
result for every nonzero sparse combination; the zero combination is
immediate.
\end{proof}

\begin{Corollary}\label{UDTsharpweighted}
Let $\Gamma\subset\bbZ^d$ be finite, $N=|\Gamma|$, and suppose
$\Phi_\Gamma=\{\phi_\bk\}_{\bk\in \Gamma}$ has a common pointwise bound
$1\le B<\infty$ on $\Omega$, with $B\in L_2(\mu)$. Assume the
lower Riesz inequality \eqref{Bessel} holds for coefficients supported
in $\Gamma$, with constant $K$. Put $c_0=\|B\|_2^{-1}$ and
$R=Kc_0^{-2}$. For each $1\le v\le N$ there are sampling points and
positive weights satisfying
\[
 m\le C R v\log(eN)\log^2(eRv),\qquad
 \sum_{i=1}^m w_i\le1,
\]
and
\[
 \frac{c_0^2}{2}\|f\|_{L_2(\mu)}^2
 \le\sum_{i=1}^m w_i|f(\xi^i)|^2
 \le\frac{3c_0^2}{2}\|f\|_{L_2(\mu)}^2,
 \qquad f\in\Sigma_v(\Phi_\Gamma).
\]
Here $C$ is absolute. The weights can be taken as
$w_i=1/(mB(\xi^i)^2)$.
\end{Corollary}
\begin{proof}
By Proposition~\ref{InP1}, the transformed system $\Psi_\Gamma$ has uniform
bound $c_0^{-1}$ and lower Riesz constant $K$ in $L_2(\nu)$. Apply
Theorem~\ref{UDTsharp} with $s=v$, choosing $m$ to satisfy
\eqref{UDsharpm} and make the success probability positive. Since
$R\ge1$, this requires only an adjustment of the absolute constant.
For $f\in\Sigma_v(\Phi_\Gamma)$, the isometry \eqref{In6} and the identity
\[
 \frac{1}{m}\sum_{i=1}^m |\alpha(\xi^i)f(\xi^i)|^2
 =c_0^{-2}\sum_{i=1}^m w_i|f(\xi^i)|^2
\]
give the assertion. Finally, $B\ge1$ implies $\sum_iw_i\le1$.
\end{proof}

We next apply this quadratic estimate to the orthonormal
basis constructed in Section~\ref{OB}.

\begin{Corollary}\label{UDRP1}
Let $\cD_N$ be an orthonormal  basis satisfying \eqref{UDRenvelope} with
$\tau=1$, and let $\rho$ be the density in \eqref{UDRsecondmeasure}.
For every $1\le v\le N$, there are sampling points $\xi^1,\ldots,\xi^m\in\Omega$ and positive weights $a_i$
such that
\begin{equation}\label{UDRpolylog}
 m\le Cv\log^2(eN)\log^2\!\bigl(ev\log(eN)\bigr)
\end{equation}
and
\[
 \frac12\|f\|_{L_2(\mu)}^2
 \le\sum_{i=1}^m a_i|f(\xi^i)|^2
 \le\frac32\|f\|_{L_2(\mu)}^2,
 \qquad f\in\Sigma_v(\cD_N).
\]
The points can be obtained as independent
samples from $\rho\,d\mu$, with $\rho$ independent of $v$, and
$a_i=1/(m\rho(\xi^i))$.
\end{Corollary}
\begin{proof}
Write $\cD_N=\{h_j\}_{j=1}^N$ and set
$g_j=h_j/\sqrt{\rho}$.
By the definition of $\rho$, the system $\{g_j\}_{j=1}^N$
is orthonormal in $L_2(\rho\,d\mu)$ and satisfies
\[
 \max_{1\le j\le N}\sup_{\bx\in\Omega}|g_j(\bx)|^2
 \le 1+\int_\Omega\max_{1\le j\le N}|h_j|^2\,d\mu
 \le C\log(eN).
\]
We may therefore apply Theorem~\ref{UDTsharp} with $K=1$
and $R=L^2\le C\log(eN)$. The required number of sampling
points is bounded by \eqref{UDRpolylog}, after rounding upward and
increasing the absolute constant to make the probability of success positive. Undoing the isometry $f\mapsto f/\sqrt\rho$ gives the displayed discretization inequality with the stated weights.
\end{proof}

For $\epsilon=1/2$, Corollary~\ref{UDRchosen} gives the smaller count
$Cv\log(eN/v)$ by
selecting evaluation points and weights. Corollary~\ref{UDRP1} uses a
single probability density for independent sampling at every sparsity.

\subsection{The bounded-system estimate for \texorpdfstring{$1\le p\le2$}{1 <= p <= 2}}
For comparison, we record the following known result; see \cite{DT,DT1}.
\begin{Theorem}[\cite{DT1}]\label{UDT1}
Let $1\le p\le2$. Suppose $\cD_N=\{g_j\}_{j=1}^N$ satisfies
$|g_j(\bx)|\le L$ on $\Omega$ and the lower Riesz inequality
\eqref{Bessel} with constant $K$. For every $1\le v\le N$ there are
$m$ points $\xi^1,\ldots,\xi^m\in\Omega$, where
\[
 m\le C v\log(eN)\log^2(2v)
       \bigl(\log(2v)+\log\log(e^2N)\bigr),
 \qquad C=C(p,L,K),
\]
such that
\begin{equation}\label{udp1}
 \frac12\|f\|_p^p\le\frac1m\sum_{i=1}^m|f(\xi^i)|^p
 \le\frac32\|f\|_p^p,\qquad f\in\Sigma_v(\cD_N).
\end{equation}
\end{Theorem}
The use of shifted logarithms allows us to include the
case $N=1$. For $N\ge2$, this form follows
from the cited estimate by increasing the constant. For $N=1$, apply
concentration to the bounded variable $|g_1|^p$, using
$\|g_1\|_p^p\ge L^{p-2}/K$. At $p=2$, Theorem~\ref{UDTsharp}
removes one logarithmic factor and records the dependence on $KL^2$.

\section{Sampling recovery}\label{SR}

For a finite dictionary $\cD_\Gamma$, we use the notation $[\cD_\Gamma]$ and
$\varrho_m^o(\bF,\cD_\Gamma,L_p(\eta))$ introduced in Section~\ref{Ex}.
Write
\[
 \sigma_v(f,\cD_\Gamma)_{L_p(\eta)}
 :=\inf_{h\in\Sigma_v(\cD_\Gamma)}\|f-h\|_{L_p(\eta)}.
\]
For a class, $\sigma_v$ denotes the supremum of this quantity over the
class. Proposition~\ref{LpRecoveryTransport} allows us to transfer the recovery
estimates for bounded systems to systems with an integrable pointwise bound.
\subsection{Sparse approximation in a finite initial dictionary}
The following approximation lemma is the $p\geq 2$ case of \cite[Theorem~1.1]{STT}. We include a proof for completeness,
recording the size bound for the finite dictionary and making
explicit that the same truncation level $J(v)$ works when the
underlying measure is replaced by a mixture of a continuous
 and an empirical measure.

\begin{Lemma}\label{SRblock}
Let $2\le p<\infty$, and suppose
$\sup_\bk\|g_\bk\|_{L_p(\eta)}\le L$ for a probability measure $\eta$.
Assume \eqref{Di4} and \eqref{sa11} hold. Put $\gamma=1/\bt-1/2$.
There are constants $C,c>0$, depending only on
$p,r,b,\bt,\th,\th'$ and $|G_1|$, such that, for every $v\ge1$, an index
$J=J(v)$ satisfies
\[
 |G_J|\le C(2v)^c,
 \qquad
 \sigma_v(\bA^{r,b}_\bt(\cD,\cG),\cD_J)_{L_p(\eta)}
 \le C L v^{-\gamma}2^{-rn(v)}n(v)^b,
 \quad \cD_J:=\{g_\bk:\bk\in G_J\}.
\]
The same $J(v)$ works for every measure satisfying the stated atom bound.
\end{Lemma}
\begin{proof}
We first prove a block estimate. If a finite coefficient vector has
$\ell_\bt$ quasi-norm at most $A$, then its expansion can be approximated
by at most $2k$ atoms with error at most $C_{p,\bt} L A k^{-\gamma}$.
For $\bt<1$, retain the $k$ largest coefficients. In decreasing order the
remaining coefficient magnitudes are bounded by $A\ell^{-1/\bt}$, so
their sum is at most $C_\bt A k^{1-1/\bt}$. Let this sum be $M$.
If $M>0$, choose an atom with probability $|a_\bk|/M$ and multiply it by
$M a_\bk/|a_\bk|$. The expectation of this random function is the tail
expansion. An average of $k$ independent copies uses at most $k$ atoms.
Symmetrization and the type $2$ estimate for $L_p$, $p\ge2$, give expected
error at most $C_p L M/\sqrt{k}$. For $\bt=1$, apply this averaging to
the entire block. Blocks of at most $2k$ entries can instead be retained
exactly.

The type $2$ estimate follows from the scalar Khintchine inequality
and Minkowski's inequality:
for $z_i\in L_p(\eta)$,
\[
 \left(\mathbb E_\varepsilon
       \left\|\sum_i\varepsilon_i z_i\right\|_p^p\right)^{1/p}
 \le C_p\left\|(\sum_i|z_i|^2)^{1/2}\right\|_p
 \le C_p(\sum_i\|z_i\|_p^2)^{1/2}.
\]
Apply this after symmetrizing the independent summands, and divide by $k$.
Hence there is a realization for which the stated block estimate holds.

Write $a_j=2^{-rj}j^b$. Put $\epsilon=r/(2\gamma)$, and choose
$c_\epsilon>0$ so that
$2c_\epsilon\sum_{\ell\ge0}2^{-\epsilon\ell}\le1/2$.
For $v\ge\max\{4|G_1|,2/c_\epsilon\}$, let
$n_0=\max\{j\ge0:|G_j|\le v/4\}$ and represent the coefficients in
$G_{n_0}$ exactly. Regularity implies
$0\le n(v)-n_0\le C_{\th'}$. Consequently
$a_{n_0+1}\asymp a_{n(v)}$, and for $\ell\ge0$,
\[
 a_{n_0+1+\ell}\le C a_{n(v)}2^{-r\ell}(1+\ell)^{|b|}.
\]
For every level with $c_\epsilon v2^{-\epsilon\ell}\ge2$, allocate
\[
 k_\ell=\lfloor c_\epsilon v2^{-\epsilon\ell}\rfloor
\]
to the block $G_{n_0+1+\ell}\setminus G_{n_0+\ell}$. The total number of
atoms used is at most $v/4+2\sum_\ell k_\ell\le3v/4$. On these levels,
$k_\ell\asymp_\epsilon v2^{-\epsilon\ell}$, so the sum of block errors is
at most
\[
 C L v^{-\gamma}a_{n(v)}
   \sum_{\ell\ge0}2^{-(r-\epsilon\gamma)\ell}(1+\ell)^{|b|}
 \le C L v^{-\gamma}a_{n(v)}.
\]
The omitted levels start at
$\ell=\epsilon^{-1}\log_2 v+O_\epsilon(1)$. Their $\ell_1$ coefficient
sum, multiplied by $L$, is bounded by
\[
 C L a_{n(v)}v^{-2\gamma}(1+\log v)^{|b|}
 \le C L a_{n(v)}v^{-\gamma}.
\]
The last retained level satisfies $J\le n_0+O(1+\log v)$, so the upper
regularity bound gives $|G_J|\le C(2v)^c$. All choices are independent of
$\eta$. For the finitely many smaller $v$, take $J=1$, use the one-term
approximant $a_{\mathbf0}g_{\mathbf0}$, and enlarge the constants. Its
error is at most $L\sum_{j\ge1}a_j$, which is finite.
\end{proof}

\subsection{Recovery for bounded systems}
The next theorem gives the recovery estimate for bounded systems.
Only a lower Riesz inequality is needed.

\begin{Theorem}\label{SRTbounded}
Let $2\le p<\infty$, and let $\eta$ be a probability measure with
$\eta(\Omega)=1$. Assume $|g_\bk(\bx)|\le L$ on $\Omega$ and
\[
 \sum_\bk|a_\bk|^2
 \le K\left\|\sum_\bk a_\bk g_\bk\right\|_{L_2(\eta)}^2
\]
for finite coefficient vectors. Under \eqref{sa11}, there are constants
$c,C>0$, depending on $p,K,L,r,b,\bt,\th,\th'$ and $|G_1|$, such that
\[
 \varrho_m^o(\bA^{r,b}_\bt(\cD,\cG),\cD_J,L_p(\eta))
 \le C v^{1-1/p-1/\bt}2^{-rn(v)}n(v)^b
\]
for $m\ge c v\log^3(2v)$, with $|G_J|\le C(2v)^c$.
\end{Theorem}
\begin{proof}
Choose $J(v)$ from Lemma~\ref{SRblock}, and put $N=|G_J|$.
Theorem~\ref{UDTsharp}, with $s=\min(2v,N)$, gives points such that
\[
 \frac{1}{2}\|d\|_{L_2(\eta)}^2
 \leq \|d\|_{L_2(\eta_m)}^2
 \le\frac{3}{2}\|d\|_{L_2(\eta)}^2,
 \qquad d\in\Sigma_s(\cD_J),
\]
where $\eta_m=m^{-1}\sum_i\delta_{\xi^i}$. Since $N\le C(2v)^c$ and
$KL^2$ is fixed, $m\ge c v\log^3(2v)$ is sufficient. Increase $c$ if
necessary to make the success probability positive.

For an $s$-sparse function $d$, the coefficient lower bound implies
\[
 \|d\|_\infty\le L\sqrt{Ks}\|d\|_{L_2(\eta)},
 \qquad
 \|d\|_{L_p(\eta)}
 \le (L\sqrt{Ks})^{1-2/p}\|d\|_{L_2(\eta)}.
\]
The second inequality follows by integrating
$|d|^p\le\|d\|_\infty^{p-2}|d|^2$.
Discretization and monotonicity of norms on the empirical probability
space therefore give
\be\label{SRsparsecompare}
 \|d\|_{L_p(\eta)}
 \le C s^{1/2-1/p}\|d\|_{L_p(\eta_m)}.
\ee

For given sample values of $f$, choose an empirical best approximant
$g_*\in\Sigma_v(\cD_J)$ in $L_p(\eta_m)$, using a fixed rule for ties.
There are finitely many supports. On each support the sample map is
injective by the lower discretization inequality, and the objective is
continuous and coercive, so a minimizer exists for every data vector.
For $h\in\Sigma_v(\cD_J)$, \eqref{SRsparsecompare} gives
\[
 \begin{split}
 \|f-g_*\|_{L_p(\eta)}
 &\le\|f-h\|_{L_p(\eta)}
       +C v^{1/2-1/p}\|h-g_*\|_{L_p(\eta_m)}\\
 &\le\|f-h\|_{L_p(\eta)}
       +2C v^{1/2-1/p}\|f-h\|_{L_p(\eta_m)}\\
 &\le C'v^{1/2-1/p}\|f-h\|_{L_p((\eta+\eta_m)/2)}.
 \end{split}
\]
The same bound $L$ holds for the dictionary elements with respect to this mixture. We now apply Lemma~\ref{SRblock}
and take the infimum over $h$. The two powers of $v$ combine as
$(1/2-1/p)+(1/2-1/\bt)=1-1/p-1/\bt$.
\end{proof}

For $p=2$, the comparison can be written with an absolute numerical
constant:
\[
 \|f-g_*\|_{L_2(\eta)}
 \le(\sqrt2+4)\sigma_v(f,\cD_J)_{L_2((\eta+\eta_m)/2)}.
\]
The reconstruction is defined by minimization over all supports.
The theorem does not provide a polynomial-time algorithm.
For results on greedy algorithms under additional hypotheses, see
\cite{VT202}.

\begin{Remark}
Under the same uniform boundedness and lower Riesz assumptions
as in \cite[Theorem~1.2]{STT}, Theorem~\ref{SRTbounded} reduces the 
sampling complexity from $m\geq cv\log^4(2v)$ to
$m\geq cv\log^3(2v)$ via Theorem~\ref{UDTsharp}, a consequence of \cite[Theorem~1.1]{BDJR}. We do not assume an upper Riesz inequality of the form 
\[
\left\|\sum_{\mathbf{k}} a_{\mathbf{k}}g_{\mathbf{k}} \right\|_{L_2(\eta)}^2 \leq C \sum_{\mathbf{k}} |a_{\mathbf{k}}|^2.
\]
\end{Remark}

\subsection{Recovery under an integrable pointwise bound}
\begin{Theorem}\label{SRT5}
Under the hypotheses of Theorem~\ref{InT1}, for every $v\ge1$ there is
$J=J(v)$ with $|G_J|\le C(2v)^c$ such that
\be\label{SR6}
 \varrho_m^o(\bA^{r,b}_\bt(\Phi,\cG),\Phi_J,L_2(\mu))
 \le C v^{1/2-1/\bt}2^{-rn(v)}n(v)^b
\ee
whenever $m\ge c v\log^3(2v)$.
\end{Theorem}
\begin{proof}
By Proposition~\ref{InP1}, the transformed system
$\Psi=\{\alpha\phi_\bk\}$ has uniform bound $c_0^{-1}=\|B\|_2$ and the same
lower Riesz constant $K$ in $L_2(\nu)$. Its coefficients are unchanged, so
\[
 \alpha\bA^{r,b}_\bt(\Phi,\cG)=\bA^{r,b}_\bt(\Psi,\cG).
\]
Apply Theorem~\ref{SRTbounded} with $p=2$, followed by
Proposition~\ref{LpRecoveryTransport}. This proves \eqref{SR6}, and also Theorem~\ref{InT1}.
For $1\le p\le2$, use $\|f\|_{L_p(\mu)}\le\|f\|_{L_2(\mu)}$.
\end{proof}

In particular, the same conclusion holds for a Riesz system satisfying
\eqref{Riesz}, with $K=R_1^{-2}$. Its upper Riesz constant is not needed.

\begin{Theorem}\label{SRTp}
Let $2<p<\infty$. Suppose \eqref{Bessel} holds and \eqref{In1} holds with
$B\in L_p(\mu)$. Under \eqref{sa11},
\[
 \varrho_m^o(\bA^{r,b}_\bt(\Phi,\cG),\Phi_J,L_p(\mu))
 \le C v^{1-1/p-1/\bt}2^{-rn(v)}n(v)^b
\]
for $m\ge c v\log^3(2v)$ and an initial dictionary satisfying
$|G_J|\le C(2v)^c$. Constants depend only on
$p,K,\|B\|_p,r,b,\bt,\th,\th'$ and $|G_1|$.
\end{Theorem}
\begin{proof}
Use $U_p$ and $\nu_p$ from \eqref{Lp3}, with $L=\|B\|_p$ and
$\alpha_p=L/B$. The atoms $\psi_\bk=U_p\phi_\bk$
have uniform bound $L$. For a finite expansion $f=\sum a_\bk\phi_\bk$,
\[
 \|U_pf\|_{L_2(\nu_p)}^2
 =L^{2-p}\int B^{p-2}|f|^2\,d\mu
 \ge L^{2-p}K^{-1}\sum_\bk|a_\bk|^2,
\]
because $B\ge1$. Thus the transformed lower Riesz constant is
$K_p=KL^{p-2}$, and the effective discretization parameter is
$K_pL^2=KL^p$. Apply Theorem~\ref{SRTbounded} to $\Psi$ with measure
$\nu_p$. Proposition~\ref{LpRecoveryTransport} transfers the recovery
error exactly. Coefficients and the structural class are preserved by
$U_p$.
\end{proof}

\begin{Corollary}\label{SRT7}
For the dyadic cube classes, under the hypotheses of Theorem~\ref{SRT5},
\[
 \varrho_m^o(\bA^r_\bt(\Phi),L_2(\mu))
 \le C v^{1/2-1/\bt-r/d},\qquad m\ge c v\log^3(2v).
\]
In particular, with $\gamma_d=r/d+1/\bt-1/2>0$,
\[
 \varrho_m^o(\bA^r_\bt(\Phi),L_2(\mu))
 \le C m^{-\gamma_d}\log^{3\gamma_d}(2m).
\]
The statements apply to orthonormal systems with a common
pointwise bound belonging to $L_2$.
\end{Corollary}
\begin{proof}
For cubes, $|G_j^c|\asymp_d2^{dj}$, whence
$2^{-rn(v)}\asymp_{r,d}v^{-r/d}$. Apply Theorem~\ref{SRT5} with $b=0$.
For large $m$ choose $v=\lfloor a m/\log^3(2m)\rfloor$ with a sufficiently
small constant $a>0$; finitely many smaller $m$ are absorbed in $C$.
\end{proof}
This is an upper estimate. We do not assert that its logarithmic factor
is optimal.

\section{Appendix: dual orthogonal systems}\label{Lp}

The change of measure also preserves the norming-functional
pairings in $L_p$. We record this fact and recall the Auerbach
construction of dual orthogonal bases. These statements describe
Banach-space structure; the sampling estimates require additional
control of the majorant and the lower $L_2$ Riesz constant.

Fix $1<p<\infty$, and use the finite systems $\Phi_\Gamma$ and $\Psi_\Gamma$
from Section~\ref{Ex}. In this appendix we write $\Phi=\Phi_\Gamma$ and
$\Psi=\Psi_\Gamma$. Recall that
\[
 \alpha_p=(c_pB)^{-1},\qquad U_pf=\alpha_pf,\qquad
 \psi_\bk=\alpha_p\phi_\bk,\qquad d\nu_p=\alpha_p^{-p}\,d\mu.
\]
By \eqref{Lp9}, $U_p$ is a surjective isometry from $L_p(\mu)$ to
$L_p(\nu_p)$.
\subsection{Coherence and norming functionals}

For a nonzero vector $g$ in a Banach space $X$, a norming functional is a
functional $F_g\in X^*$ with $\|F_g\|_{X^*}=1$ and
$F_g(g)=\|g\|_X$. If $X$ is smooth, this functional is unique. In particular, norming functionals are unique in $L_p$ when $1<p<\infty$.
For a normalized system $\cD\subset X$, define its coherence by
\[
 M(\cD,X):=\sup_{g\ne h;\ g,h\in\cD}|F_g(h)|.
\]
For general Banach spaces one also takes the supremum over all norming
functionals of $g$; see \cite[p.~381]{VT107}. Here we work in smooth spaces,
so this additional supremum is unnecessary. Clearly $0\le M(\cD,X)\le1$.
For an ordered normalized system $\cD_N=\{g^j\}_{j=1}^N$, the upper and
lower semi-incoherence parameters introduced in \cite{VT206} are
\[
 M^u(\cD_N,X):=\sup_{1\le i<j\le N}|F_{g^j}(g^i)|,
 \qquad
 M^\ell(\cD_N,X):=\sup_{1\le i<j\le N}|F_{g^i}(g^j)|.
\]
Reversing the order interchanges these two parameters. The upper parameter
also makes sense for an infinite ordered system. A supremum over an empty
set of pairings is understood to be zero.

Remove any zero elements from $\Phi$ and define
\[
 \widehat\Phi:=\left\{\frac{\phi_\bk}{\|\phi_\bk\|_{L_p(\mu)}}\right\}_{\bk\in \Gamma},
 \qquad
 \widehat\Psi:=\left\{\frac{\psi_\bk}{\|\psi_\bk\|_{L_p(\nu_p)}}\right\}_{\bk\in \Gamma}.
\]
This normalization is needed to define the coherence parameters. The corresponding denominators are equal by \eqref{Lp9}.

\begin{Proposition}\label{LpP1}
Let $1<p<\infty$. Under the assumptions of \eqref{Lp1}--\eqref{Lp3},
\be\label{Lp4}
 M(\widehat\Phi,L_p(\mu))=M(\widehat\Psi,L_p(\nu_p)).
\ee
In particular, the equality holds without hats if $\Phi$ is normalized in
$L_p(\mu)$.
\end{Proposition}

\begin{Proposition}\label{LpP2}
Let $1<p<\infty$, and order $\Phi$ and $\Psi$ in the same way. Then
\be\label{Lp5}
 \begin{aligned}
 M^u(\widehat\Phi,L_p(\mu))&=M^u(\widehat\Psi,L_p(\nu_p)),\\
 M^\ell(\widehat\Phi,L_p(\mu))&=M^\ell(\widehat\Psi,L_p(\nu_p)).
 \end{aligned}
\ee
\end{Proposition}

\begin{proof}
For a nonzero $f\in L_p(\mu)$, its norming functional is
\be\label{Lp6}
 F_f^{(p,\mu)}(g)
 =\|f\|_{L_p(\mu)}^{1-p}
   \int_\Omega |f|^{p-2}\overline f g\,d\mu,
\ee
where $|f|^{p-2}\overline f$ is taken to be zero at $f=0$.
Since $\alpha_p$ is positive, \eqref{Lp9} and $\alpha_p^pc_p^pB^p=1$ give
\be\label{Lp7}
 F_{U_pf}^{(p,\nu_p)}(U_pg)=F_f^{(p,\mu)}(g).
\ee
Applying this identity to every pair of normalized elements proves both
propositions.
\end{proof}

\subsection{Construction of dual orthogonal systems}

\begin{Definition}\label{LpD1}
Let $X$ be a smooth Banach space. A normalized finite or infinite system
$\cD=\{g^j\}$ is \emph{dual orthogonal} if
$F_{g^j}(g^k)=0$ for every $j\ne k$. We write $\cD\in DOS(X)$.
\end{Definition}

\begin{Definition}\label{LpD2}
A normalized ordered system $\cD=\{g^j\}$ in a smooth Banach space $X$ is
\emph{upper dual orthogonal} if $F_{g^j}(g^k)=0$ whenever $j>k$.
We write $\cD\in UDOS(X)$.
\end{Definition}

Identity \eqref{Lp7} shows that $U_p$ preserves both kinds of dual
orthogonality. The same conclusion follows for any surjective linear
isometry between smooth Banach spaces, since its adjoint transports the
unique norming functionals.

The following construction gives upper dual orthogonality while preserving
the successive spans. Let $x^1,\ldots,x^N$ be linearly independent in $X$.
Set $g^1=x^1/\|x^1\|_X$. For $j\ge2$, let $y^j$ be a best approximant to
$x^j$ from the finite-dimensional space
$E_{j-1}:=\operatorname{span}\{x^1,\ldots,x^{j-1}\}$, and put
\[
 g^j:=\frac{x^j-y^j}{\|x^j-y^j\|_X}.
\]
A best approximant exists by finite-dimensional compactness, and the
residual is nonzero by linear independence. The first-order optimality
condition for the norm at this residual gives
$F_{g^j}(z)=0$ for every $z\in E_{j-1}$. In the complex case, apply the
real directional derivative to $z$ and $iz$. Inductively,
$\operatorname{span}\{g^1,\ldots,g^{j-1}\}=E_{j-1}$, so the resulting
system belongs to $UDOS(X)$.

To obtain dual orthogonality, we optimize all the vectors simultaneously.

\begin{Theorem}[Auerbach construction]\label{LpAuerbach}
Let $E$ be an $N$-dimensional subspace of a smooth real or complex Banach
space $X$. Then there is a basis $g^1,\ldots,g^N$ of $E$ such that
\[
 \|g^i\|_X=1,\qquad F_{g^i}(g^j)=\delta_{ij},
 \qquad 1\le i,j\le N.
\]
In particular, every such subspace admits a dual orthogonal basis.
\end{Theorem}
\begin{proof}
Fix linear coordinates on $E$, and maximize
$|\det(g^1,\ldots,g^N)|$ over the $N$-fold product of its closed unit
ball. Compactness gives a maximizer, and its determinant is nonzero because
$E$ has a basis. Every maximizing column has norm one: otherwise scaling
that column to norm one would increase the determinant.

Let $g_i^*$ be the coordinate functional of this maximizing basis.
Replacing column $i$ by any $x$ in the unit ball of $E$ multiplies the
determinant by $g_i^*(x)$. Maximality therefore gives $|g_i^*(x)|\le1$.
Since $g_i^*(g^i)=1$, we have $\|g_i^*\|_{E^*}=1$. The Hahn--Banach
theorem extends $g_i^*$ to a norm-one functional on $X$ taking the value one at $g^i$.
Smoothness makes this extension the unique norming functional $F_{g^i}$.
The coordinate identities $g_i^*(g^j)=\delta_{ij}$ complete the proof.
\end{proof}
\begin{Remark}
The constructed basis and its coordinate functionals form an
Auerbach basis of $E$; that is,
\[
\|g^i\|_E=\|g_i^*\|_{E^*}=1,
\qquad
g_i^*(g^j)=\delta_{ij},
\qquad 1\le i,j\le N.
\]
For $E\subset L_p(\mu)$, $1<p<\infty$, the theorem yields the equations
\[
 \|g^i\|_{L_p(\mu)}=1,\qquad
 \int_\Omega |g^i|^{p-2}\overline{g^i}g^j\,d\mu=\delta_{ij}.
\]
Thus the finite-dimensional construction can be formulated as a
maximum-determinant problem over linear combinations of any given basis of
$E$.

\end{Remark}

A prescribed initial vector cannot always be retained in this construction.
For example, in $\ell_4^2$ take the direction $(1,2)$. The condition
$F_{(1,2)}(g)=0$ forces $g$ to have direction $(8,-1)$. The reverse pairing
is nonzero, because
\[
 8^3\cdot1+(-1)^3\cdot2=510\ne0.
\]
Normalization does not affect these zero conditions. Consequently, a successive residual construction does not in general give dual orthogonality while retaining an arbitrary first vector.

There are also explicit infinite examples in $L_p$, $1<p<\infty$.
Normalized functions with pairwise disjoint supports form a dual orthogonal
system. So do orthogonal functions of constant modulus one, including the
Rademacher and Walsh systems on $[0,1]$ and the distinct characters of a
compact abelian group with normalized Haar measure. The $L_p$-normalized
dyadic Haar functions on $[0,1]$, together with the constant function, give
another example. Indeed, for a Haar function $h_I$ supported on $I$,
\[
 |h_I|^{p-2}h_I=|I|^{-(p-2)/p}h_I,
\]
and distinct dyadic Haar functions have zero ordinary integral pairing.
The infinite Haar system has no common finite pointwise bound: the
amplitude $|I|^{-1/p}$ grows along the dyadic intervals containing almost
every fixed point.

More generally, Day's theorem gives an infinite Auerbach system in every
infinite-dimensional Banach space; see \cite[Remark~3.5]{HKR}. In a smooth
space the associated norm-one biorthogonal functionals are the unique
norming functionals, so this yields an infinite dual orthogonal system.
This existence statement does not assert that the system is a basis of the
whole space, or that it has a prescribed order, regularity, or common integrable pointwise bound.

\subsection{Consequences and limitations for sampling}

For a normalized dual orthogonal system and a finite scalar sequence,
\[
 \max_j|a_j|
 \le\left\|\sum_j a_jg^j\right\|_X
 \le\sum_j|a_j|.
\]
The first inequality follows by applying $F_{g^i}$, and the second is the
triangle inequality. These bounds do not imply the lower $\ell_2$
estimate needed in the sampling results of the preceding sections. On $N$
equally weighted atoms, for example, the functions
$g^j=N^{1/p}\mathbf1_{\{j\}}$ are normalized and dual orthogonal in $L_p$,
but
\[
 \left\|\sum_{j=1}^N a_jg^j\right\|_{L_2}^2
 =N^{2/p-1}\sum_{j=1}^N|a_j|^2.
\]
For $p>2$ their lower Riesz constant therefore grows as $N^{1-2/p}$.
Likewise, an orthogonal change of coefficients preserves $L_2$
orthogonality but need not preserve dual orthogonality in $L_p$ for
$p\ne2$. Thus the constructions in Sections~\ref{UDR} and~\ref{OB} require separate
control of the sampling estimates and of the coefficient classes.

\section*{Acknowledgment}
The first and second authors were partially supported by
NSERC of Canada Discovery Grant RGPIN~2026-05667. The third author’s research (Sections 4 and 5) was carried out with
the financial support of the Ministry of Science and Higher Education
of the Russian Federation in the framework of a scientific project
under agreement No. 075-15-2025-013.

\Addresses

\end{document}